\documentclass[11pt,a4paper,reqno]{amsart}
\usepackage{graphicx}
\usepackage{subfigure}
\usepackage{comment}
\usepackage{tikz, pgfplots}
\usepackage{hyperref}
\usepackage{epstopdf}

\hypersetup{
colorlinks=true,
linkcolor=blue,
anchorcolor=blue,
citecolor=blue
}

\usepackage{amssymb}
\usepackage{amsfonts}
\usepackage{amsmath}
\usepackage{mathrsfs}
\usepackage{color}

\newtheorem{theorem}{Theorem}[section]
\newtheorem{lemma}[theorem]{Lemma}
\newtheorem{proposition}[theorem]{Proposition}
\newtheorem{corollary}[theorem]{Corollary}

\theoremstyle{definition}

\newtheorem{remark}[theorem]{Remark}

\allowdisplaybreaks

\begin{document}
\title[quasilinear elliptic equations]{A priori estimates  and Liouville type results for quasilinear elliptic equations\\ involving gradient terms}
\author[Filippucci]{Roberta Filippucci}
\address{Dipartimento di Matematica e Informatica, Universit\`{a} degli Studi di Perugia, Via Vanvitelli 1, 06123 Perugia, Italy}
\email{roberta.filippucci@unipg.it}
\thanks{Filippucci is a member of the {\em Gruppo Nazionale per
l'Analisi Ma\-te\-ma\-ti\-ca, la Probabilit\`a e le loro Applicazioni}
(GNAMPA) of the {\em Istituto Nazionale di Alta Matematica} (INdAM)
 and was  partly supported by  {\em Fondo Ricerca di
Base di Ateneo Esercizio} 2017-19 of the University of Perugia
"Problemi con non linearit\`a dipendenti dal gradiente" and by  GNAMPA-INdAM Project 2022 "Equazioni differenziali alle derivate parziali in fenomeni non lineari"
(CUP-E55F22000270001)}
\author[Sun]{Yuhua Sun}
\address{School of Mathematical Sciences and LPMC, Nankai University, 300071 Tianjin, P. R. China}
\email{sunyuhua@nankai.edu.cn}
\thanks{ Sun was supported by the National Natural Science Foundation of China (No.11501303).}

\author[Zheng]{Yadong Zheng}
\address{School of Mathematical Sciences and LPMC, Nankai University, 300071 Tianjin, P. R. China}
\email{yadongzheng2017@sina.com}



\begin{abstract}
In this article we study local and global properties of positive solutions of $-\Delta_mu=|u|^{p-1}u+M\left|\nabla u\right|^q$ in a domain $\Omega$ of $\mathbb R^N$, with $m>1$, $p,q>0$ and $M\in\mathbb R$. Following some ideas used in \cite{BV,Vron1}, and by using a direct Bernstein method combined with Keller-Osserman's estimate, we obtain several a priori estimates as well as Liouville type theorems. Moreover, we prove a local Harnack inequality with the help of Serrin's classical results.
\end{abstract}

\maketitle

\section{Introduction}
In this paper, we aim to investigate local and global properties of positive solutions to the following equation
\begin{equation}\label{1-1}
-\Delta_mu=|u|^{p-1}u+M\left|\nabla u\right|^q\quad{\rm in}\;\Omega,
\end{equation}
where $m>1$, $\Delta_mu={\rm div}\left(\left|\nabla u\right|^{m-2}\nabla u\right)$, $p,q>0$, $M\in\mathbb R$
and $\Omega\subset\mathbb R^N$ ($N\geq1$) is a domain bounded or not and containing $0$.

If $M=0$, then  \eqref{1-1} reduces to the generalized Lane-Emden equation
\begin{equation}\label{0-1}
-\Delta_mu=|u|^{p-1}u\quad{\rm in}\;\Omega,
\end{equation}
which has been widely studied in the literature \cite{Bidaut,Bidaut1,Caffarelli,Chen,Damascelli,Gidas,Guedda,Ni,NL,Sciunzi,Serrin,SZ,Vtois}, both when $\Omega$ is bounded
and when $\Omega$ is unbounded.
Especially, in the semilinear case $m=2$, one of the celebrated results is given by Gidas and Spruck \cite{Gidas}:
if $N>2$ and $p\in\left[1,\frac{N+2}{N-2}\right)$, then any nonnegative solution of (\ref{0-1}) in $\mathbb R^N$ is identically zero and the result is sharp.
Very surprisingly in Gidas-Spruck's result, there is no a priori information assumption on the behavior of the solutions at infinity. Additional results for the semilinear case, but with a nonlinearity similar to that  in \eqref{1-1}
can be found in \cite{DuGheRad}  and \cite{GheRad}.

For the case of $m>1$,  radially symmetric positive solutions
 were studied by Ni and Serrin \cite{Ni,NL,Serrin1}, and  further results in this direction were obtained by Guedda and V$\acute{\rm e}$ron \cite{Guedda} and Bidaut-V$\acute{\rm e}$ron \cite{Bidaut}.

When one studies the so called Liouville property of (\ref{0-1}), namely whether all  positive $C^1$ solutions of (\ref{0-1}) in $\mathbb R^N$ are constant,
 two critical exponents appear
\begin{equation}\label{0-2}
m_*=\frac{N(m-1)}{N-m},\quad m^*=\frac{N(m-1)+m}{N-m},
\end{equation}
when $N>m$, known as the Serrin exponent and the Sobolev exponent, respectively.  It is well known that the first is optimal for the Liouville property for the inequality
\begin{equation}\label{0-4}
-\Delta_mu\geq|u|^{p-1}u,\quad{\rm in}\;\mathbb R^N,
\end{equation} while the second is optimal  for the corresponding equality.  Indeed, Mitidieri and Pohozaev \cite{Mitidieri} first proved that if $N>m$ and $p\in(0, m_*]$, or $N\leq m$ and $p\in(0, \infty)$, then any nonnegative solution to \eqref{0-4}  is zero. On the other hand, if $N>m$ and $p\in(m_*, \infty)$, then \eqref{0-4} possesses the following bounded positive solution
$$u(x)=C\left(1+|x|^{\frac{m}{m-1}}\right)^{-\frac{m-1}{p-m+1}},$$
for some $C>0$, see \cite[Remark 4]{Mitidieri} or \cite{SZ}.
For equation  \eqref{0-1}  in $\mathbb R^N$, we refer to the marvellous paper by Serrin and Zou \cite{SZ} (cfr. Corollary II), where also  nonexistence in the case $N<m$ and $p\in(0,\infty)$ was solved. Of course, if $M\geq0$, every positive solution of (\ref{1-1})  is also a positive solution of the inequality \eqref{0-4}.

If we consider the critical case of  (\ref{0-1}), that is when $p=m^*$, and we restrict our attention to solutions belonging to  the space $\mathcal D^{1,m}(\mathbb R^N):=\left\{u\in L^{m^*+1}(\mathbb R^N):\int_{\mathbb R^N}|\nabla u|^m<\infty\right\}$,
then Damascelli et al. in \cite{Damascelli}, for $1<m<2$,  Sciunzi in \cite{Sciunzi}, for $m>2$ , and  V\`etois in \cite{Vtois}, for  $m>1$, showed that all positive solutions  are radial and have the following form
$$
u(x)=U_{\lambda,x_0}(x):=\left[\frac{\lambda^{\frac{1}{m-1}}N^{\frac{1}{m}}\left(\frac{N-m}{m-1}\right)^{\frac{m-1}{m}
}}{\lambda^{\frac{m}{m-1}}+|x-x_0|^{\frac{m}{m-1}}}\right]^{\frac{N-m}{m}},\quad\lambda>0,\;x_0\in\mathbb R^N.
$$

Moving  to exterior domains,  Bidaut-V$\acute{\rm e}$ron \cite{Bidaut} proved that any nonnegative solution of (\ref{0-1}) is zero provided that $N>m$ and $p\in(m-1, m_*]$, or $N=m$ and $p\in(m-1, \infty)$, while  Bidaut-V$\acute{\rm e}$ron and Pohozaev \cite{Pohozaev} showed that (\ref{0-4}) admits only the trivial solution $u\equiv0$ whenever $N>m$ and $p\in(0,m_*]$, or $N=m$ and $p\in(0,\infty)$.

For the case with gradient terms, we first recall the Hamilton-Jacobi equation
\begin{equation}\label{0-5}
-\Delta_mu=\left|\nabla u\right|^q\quad{\rm in}\;\Omega.
\end{equation}
The Liouville property of \eqref{0-5} was studied by Lions in \cite{Lions} for $m=2$, who proved that any $C^2$ solution
to \eqref{0-5} with $q>1$ in $\mathbb R^N$ has to be a constant by using the Bernstein technique. Bidaut-V$\acute{\rm e}$ron, Garcia-Huidobro and V$\acute{\rm e}$ron \cite{Bidaut2} proved that any $C^1$ solution $u$ of (\ref{0-5}) in an
arbitrary domain $\Omega$ of $\mathbb R^N$ with $N\geq m>1$ and $q>m-1$ satisfies
\begin{equation}\label{0-6}
\left|\nabla u(x)\right|\leq c_{N,m,q}
\left({\rm dist}\left(x,\partial\Omega\right)\right)^{-
\frac{1}{q-m+1}}
\end{equation}
for all $x\in\Omega$.  Estimates of this type, not only for the gradient but also for the solutions are called by Serrin and Zou "universal a priori estimates", because they are independent of the solutions and do not need any boundary conditions. In particular, they   produce as  a direct
corollary the Liouville property, since ${\rm dist}\left(x,\partial\Omega\right)$ can be chosen arbitrarily large when the solution is defined on all $\mathbb R^N$. For a detailed discussion in this direction we refer to the paper by  Polacik, Quitter and Souplet studied in \cite{PQS} where new connections between Liouville type theorems and universal estimates were developed.
Here ``any solution" means there is no any sign condition on the solution. Estimates of the gradient for more general problems can be found in \cite{Leonori}.

For the generalized case of \eqref{0-5} given by
\begin{equation}\label{0-7}
-\Delta_mu=u^p\left|\nabla u\right|^q\quad{\rm in}\;\Omega,
\end{equation}
in \cite{Ver} Bidaut-V$\acute{\rm e}$ron, Garcia-Huidobro and V$\acute{\rm e}$ron focused on positive solutions of (\ref{0-7}) for $m=2$, $p\geq0$ and $0\leq q<2$. By using the pointwise Bernstein method and the integral Bernstein method, they determined various regions of $(p,q)$ for which the Liouville property holds.
Filippucci, Pucci and Souplet \cite{Filippucci} solved the case of $m=2$, $p>0$ and $q>2$, and they proved that any positive bounded classical solution of (\ref{0-7}) in $\mathbb R^N$ is identically equal to a constant.
Bidaut-V$\acute{\rm e}$ron \cite{Bidaut3} obtained the same Liouville-type results for (\ref{0-7}) in the case $N> m>1$, $p\geq0$ and $q\geq m$ without the assumption of boundedness on the solution.
Recently, the Liouville property of (\ref{0-7}) in $\mathbb R^N$ for $N\geq1$, $m>1$, $p\geq0$ and $0\leq q<m$ was studied by Chang, Hu and Zhang \cite{Chang}.
For the case of  radial solutions of the coercive vectorial version of  (\ref{0-7}) in $\mathbb R^N$ we refer to \cite{GGS}.

If we consider the inequality version of (\ref{0-7})
\begin{equation}\label{0-77}
-\Delta_mu\geq u^p\left|\nabla u\right|^q\quad{\rm in}\;\Omega,
\end{equation}
it was  proved in \cite{Ver}   for the case $m=2$ that any positive solution of (\ref{0-77}) in $\mathbb R^N$ must be constant if
$N> 2$, $p\geq0$, $q\geq0$ and $$p(N-2)+q(N-1)<N.$$ The generalization of the above results to the case $m\neq2$, even in the vectorial case can be found in \cite{RF, RF1, RF2, Mitidieri1}.

Recently, Sun, Xiao and Xu \cite{Sun} dealt with (\ref{0-77}) when $\Omega$ is a geodesically complete noncompact Riemannian manifold,  and they obtained the nonexistence and existence of positive solutions to (\ref{0-77}) in the range $m>1$ and $(p,q)\in\mathbb R^2$ via the volume growth of geodesical ball.

The most important motivation of the present study is to extend the results obtained for the semilinear equation
\begin{equation}\label{0-8}
-\Delta u=|u|^{p-1}u+M\left|\nabla u\right|^q\quad{\rm in}\;\Omega,
\end{equation}
by Bidaut-V$\acute{\rm e}$ron, Garcia-Huidobro and V$\acute{\rm e}$ron, see \cite{BV,Vron1}. By using a delicate combination of refined Bernstein techniques and Keller-Osserman estimate, they obtained a series of a priori estimates for any positive solution of (\ref{0-8}) in  arbitrary domain $\Omega$ of $\mathbb R^N$ in the
case  $p>1$, $q\geq\frac{2p}{p+1}$ and $M>0$ (\cite[Theorems A, C, D]{BV}). In particular the nonexistence of positive solutions of (\ref{0-8}) in $\mathbb R^N$ was obtained for the following cases:
\begin{enumerate}
\item[(i)]{$N\geq1$, $p>1$, $1<q<\frac{2p}{p+1}$, $M>0$ ;}
\item[(ii)]{$N\geq1$, $p>1$, $q=\frac{2p}{p+1}$, $M>\left(\frac{p-1}{p+1}\right)^{\frac{p-1}{p+1}}\left(\frac{N(p+1)^2}{4p}\right)^{\frac{p}{p+1}}$;}
\item[(iii)]{$N\geq2$, $1<p<\frac{N+3}{N-1}$, $1<q<\frac{N+2}{N}$, $M>0$;}
\item[(iv)]{$N\geq3$, $1<p<\frac{N+2}{N-2}$, $q=\frac{2p}{p+1}$, $|M|\leq\epsilon_0$,}
\end{enumerate}
where $\epsilon_0$ is a positive constant given in \cite[Theorem E]{BV}.
They also considered the existence and nonexistence of ``large solutions",  namely those solutions $u(x)\to\infty$ as ${\rm dist}\left(x,\partial\Omega\right)\to0$, and radial solutions of (\ref{0-8}).

In this paper, we follow the idea used in \cite{BV,Vron1}, based on the Bernstein method, to derive various a priori estimates concerning $\nabla u$ for  positive solutions of \eqref{1-1}
 in the cases $q$ is less, greater  or equal to $\frac{mp}{p+1}$,  and consequently we obtain  Liouville  type  theorems.

 Our first result is devoted to the case $q>\frac{mp}{p+1}$.
\begin{theorem}\label{thm1}
{\rm Let $\Omega\subset\mathbb R^N$, $p>\max\left\{m-1,1\right\}$ and $q>\frac{mp}{p+1}$.

Then for any $M>0$, there exists a positive constant $c_{N,m,p,q}$ such that any positive solution of (\ref{1-1}) in $\Omega$ satisfies
\begin{equation}\label{1-2}
\left|\nabla u(x)\right|\leq c_{N,m,p,q}\left(M^{-\frac{p+1}{(p+1)q-mp}}
+\left(M{\rm dist}\left(x,\partial\Omega\right)\right)^{-\frac{1}{q-m+1}}\right)\end{equation}
for all  $x\in\Omega$.
Especially, any positive solution of (\ref{1-1}) in $\mathbb R^N$ has at most a linear growth at infinity
\begin{equation}\label{1-3}
\left|\nabla u(x)\right|\leq c_{N,m,p,q}M^{-\frac{p+1}{(p+1)q-mp}},\quad x\in\mathbb R^N.
\end{equation}
}
\end{theorem}

\medskip While in the case $q<\frac{mp}{p+1}$, we obtain a nonexistence result.
\begin{theorem}\label{thm2}
{\rm Let $p>\max\left\{m-1,1\right\}$ and $\max\left\{m-1,\frac{m}{2}\right\}<q<\frac{mp}{p+1}$.

Then for any $M>0$, there exists a positive constant $c_{N,m,p,q}$ such that  (\ref{1-1})  does not admit positive solutions in $\mathbb R^N$ satisfying
\begin{equation}\label{1-4}
u(x)\leq c_{N,m,p,q}M^{\frac{m}{mp-(p+1)q}},\quad x\in\mathbb R^N.
\end{equation}
}
\end{theorem}

\begin{remark}
Here the condition $q>\frac{m}{2}$ is necessary from Young's inequality, otherwise Theorem \ref{thm2} is not valid
any more.
\end{remark}

For the case $q=\frac{mp}{p+1}$ and $M$ large enough, we have the following nonexistence result in $\mathbb R^N$.

\begin{theorem}\label{thm3}
{\rm Let $\Omega\subset\mathbb R^N$, $p>\max\left\{m-1,1\right\}$ and $q=\frac{mp}{p+1}$.

Then for any
\begin{equation}\label{1-5}
M>\frac{\sqrt{N}(p+1)}{(4p)^{\frac{p}{p+1}}}\left(\frac{p-1}{\sqrt{a}}\right)^{\frac{p-1}{p+1}},
\end{equation}
where $0<a\leq\frac{1}{N}$, there exists a positive constant $c_{N,M,a,m,p,q}$ such that any positive solution of (\ref{1-1}) in $\Omega$ satisfies
\begin{equation}\label{1-5-5}
\left|\nabla u(x)\right|\leq c_{N,M,a,m,p,q}
\left({\rm dist}\left(x,\partial\Omega\right)\right)^{-
\frac{p+1}{p-m+1}}\end{equation}
for all $ x\in\Omega$. Consequently, (\ref{1-1})  does not admit positive solutions in $\mathbb R^N$ }.
\end{theorem}

When $M$ is allowed to be negative, we derive a nonexistence result for supersolutions of (\ref{1-1}) in an exterior domain.

\begin{theorem}\label{thm4}
{\rm Let $p>m-1$ if $N= m$
or
$m-1<p<\frac{N(m-1)}{N-m}$ if $N> m$, $q=\frac{mp}{p+1}$ and $M>-\mu^\ast(N)$
where
\begin{equation}\label{1-6}
\mu^\ast(N):=(p+1)\left(\frac{N(m-1)-p(N-m)}{mp}\right)^{\frac{p}{p+1}}.
\end{equation}
Then there exist no nontrivial nonnegative supersolutions of (\ref{1-1}) in  $\mathbb R^N\setminus \overline{B}_R $ for any $R>0$.}
\end{theorem}

Concerning  large solutions,  we prove the following.

\begin{theorem}\label{thm5}
{\rm Let $\Omega$ be a open domain with Lipschitz boundary,  $p>m-1$ and $q=\frac{mp}{p+1}$.
If $M\geq-\mu^\ast(m)$, then there exists no positive supersolution of (\ref{1-1}) in $\Omega$ satisfying
\begin{equation}\label{1-7}
\lim_{{\rm dist}\left(x,\partial\Omega\right)\rightarrow0}u(x)=\infty.
\end{equation}
}
\end{theorem}

Inspired by \cite[Theorem A]{Ver}, we derive an a priori estimate for positive solution $u$ of (\ref{1-1}) in a neighborhood of $0$ as follows. 
The proof relies on Serrin's classical Harnack inequality \cite[Theorem 5]{Serrin} and the fact that  every  radial solution  $u(|x|)$ of
(\ref{1-1}) is $m$-superharmonic when $M\geq0$.

\begin{theorem}\label{thm0}
{\rm Let $\Omega\subset\mathbb R^N$ ($N\geq2$) be a domain containing $0$. Assume $1<m<N$, $m-1<p<\frac{N(m-1)}{N-m}$, $m-1<q<\frac{N(m-1)}{N-1}$
and $M\geq0$. If $u\in C^2\left(\Omega\backslash\{0\}\right)$ is a positive solution of (\ref{1-1}) in $\Omega\backslash\{0\}$, then
\begin{equation}\label{00}
u(x)+|x||\nabla u(x)|\leq c|x|^{\frac{m-N}{m-1}}
\end{equation}
holds in a neighborhood of $0$ for some $c>0$.
}
\end{theorem}

\begin{remark}
Under the assumptions on $N,m,p,q$ and $M$ of Theorem \ref{thm0}, we obtain a local Harnack inequality for positive solution $u$ of (\ref{1-1}), namely
\begin{align}\label{000-1}
\max_{|x|=r}u(x)\leq K\min_{|x|=r}u(x),\quad r\in\left(0,1/2\right],
\end{align}
for some $K>0$. The Harnack inequality for more general model
\begin{align}\label{000-2}
|u|^{p-1}u-M\left|\nabla u\right|^q\leq-\Delta_mu\leq c_0|u|^{p-1}u+M\left|\nabla u\right|^q,
\end{align}
where $c_0\geq1$ and $M>0$, was obtained first by Ruiz \cite{R} in the range $m-1<p<\frac{N(m-1)}{N-m}$ and $m-1<q<\frac{mp}{p+1}$. Note here $\frac{mp}{p+1}<\frac{N(m-1)}{N-1}$ 
always holds if $p$ satisfies the assumption of Theorem \ref{thm0}. 
\end{remark}

The final result is a Liouville-type theorem for positive solution of \eqref{1-1} with a less restrictive assumption on $M$ but a more restrictive assumption on $p$ compared with Theorem \ref{thm3}.
 Actually, as emphasized before  \cite[Theorem B]{Ver}, the direct Bernstein method allows to obtain pointwise estimates of the gradient without any integration. In particular, in the next result, using
 cumbersome algebraic manipulations and a rather demanding application of Young's inequality, we obtain an a priori estimate for the norm of the gradient of a power of a positive solution, in the spirit of  \cite[Theorem B]{Ver} devoted to elliptic inequality  of the Laplacian type with a superlinear absorption term.

\begin{theorem}\label{thm6}
{\rm Let $\Omega\subset\mathbb R^N$ ($N\geq2$). Assume $m-1<p<\frac{(N+3)(m-1)}{N-1}$ and $m-1<q<\frac{(N+2)(m-1)}{N}$. Then for any $M>0$, there exist positive constants $d$ and
$c_{N,m,p,q}$ such that any positive solution of (\ref{1-1}) in $\Omega$ satisfies
\begin{equation}\label{1-8}
\left|\nabla u^d(x)\right|\leq c_{N,m,p,q}
\left({\rm dist}\left(x,\partial\Omega\right)\right)^{-1-\frac{md}{p-m+1}},\quad x\in\Omega.
\end{equation}
In particular, there exists no nontrivial nonnegative solution of (\ref{1-1}) in $\mathbb R^N$.
}
\end{theorem}

As a consequence of (\ref{1-8}) the following holds, we have
\begin{corollary}
\rm{
Let $\Omega$ be a smooth domain in $\mathbb R^N$ ($N\geq2$) with a bounded boundary, and under the assumptions of Theorem \ref{thm6}. If $u$ is a positive solution of (\ref{1-1}) in $\Omega$, then there exists a positive constant $d_0$ depending on $\Omega$ and $c_{N,m,p,q}>0$ such that
\begin{equation}\label{6-1}
u(x)\leq c\left(\left({\rm dist}\left(x,\partial\Omega\right)\right)^{-
\frac{m}{p-m+1}}+\max_{{\rm dist}(z,\partial\Omega)=d_0} u(z)\right),\quad x\in\Omega.
\end{equation}}
\end{corollary}

\textbf{Notations.} In the above and below, the letters $C,C',C_0,C_1,c_0,c_1$... denote positive constants whose values are unimportant and may vary at different occurrences, and $C_{x,\cdots,z}$ or $C(x,\cdots, z)$ denotes the positive constant whose value relies on the choices of $x,\cdots, z$.

\section{Proof of Theorems \ref{thm1}, \ref{thm2} and \ref{thm3}}
We begin with the following lemma which plays a key role in our proofs.
\begin{lemma}\label{lem1}
\rm{ Let $\Omega\subset\mathbb R^N$, $N\geq1$ and $m>1$. Assume that $v$ is a $C^1$ function in $\Omega$  such that $|\nabla v|>0$, and let $w$ be a  continuous and nonnegative  function in $\Omega$ with  $w\in C^{ 2}(\mathcal{W}_+)$, 
where $\mathcal{W}_+=\left\{x\in\Omega:w(x)>0\right\}$. Define the operator
$$w\rightarrow\mathscr{A}_v(w):=-\Delta w-(m-2)\frac{\left<D^2w\nabla v,\nabla v\right>}{|\nabla v|^2}.$$
If $w$ satisfies, for some $\xi>1$ and a real number $c_0$,
$$\mathscr A_v(w)+w^\xi\leq c_0\frac{\left|\nabla w\right|^2}{w}$$
on each connected component of $\mathcal{W}_+$, then
$$w(x)\leq c_{N,\xi,c_0}\left({\rm dist}\left(x,\partial\Omega\right)\right)^{-\frac{2}{\xi-1}},\quad\forall x\in\Omega.$$
In particular, $w\equiv0$ if $\Omega=\mathbb R^N$.
}
\end{lemma}
\rm \begin{proof}
This proof is a combination of  \cite[Proposition 2.1]{Bidaut2},  and that of   \cite[Lemma 3.1]{Bidaut3} in the special case $\beta(x)=0$.  In particular, the operator $\mathscr{A}_v(w)$ was first introduced in  \cite[Proposition 2.1]{Bidaut2}.
\end{proof}
  The next lemma is the extension of  formula (2.6)  in \cite{BV}, the new formula, valid for every $m>1$,  is rather tricky and requires cumbersome calculations since  we have to take into account several terms appearing when $m\neq2$. 
\begin{lemma}\label{lem3}
	\rm{ Assume that $v$ is  a nonnegative 
		 $C^{3+\alpha}$ function  in $\Omega$ for some $\alpha\in(0,1)$. Let $z=|\nabla v|^2$, then we have
\begin{align}\label{L-1}
	&\quad\frac{1}{2}\mathscr A_v(z)+\frac{1}{N}z^{2-m}\left(\Delta_mv\right)^2
	+z^{1-\frac{m}{2}}\left<\nabla\Delta_mv,\nabla v\right>
	\nonumber\\&
	\leq \frac{(N+2)(m-2)}{2N}z^{-\frac{m}{2}}\Delta_mv\left<\nabla z,\nabla v\right>+\frac{m-2}{4}\frac{|\nabla z|^2}{z}
	\nonumber\\&\quad
	-\frac{(2N+m-2)(m-2)}{4N}\frac{\left<\nabla z,\nabla v\right>^2}{z^2},\quad\mbox{on $\{z>0\}$}.
\end{align}
}
\end{lemma}
\rm \begin{proof}
	Using $z=|\nabla v|^2$, $\nabla z=2D^2v\nabla v$, and
	\begin{align}
		\Delta_mv=|\nabla v|^{m-2}\Delta v+(m-2)|\nabla v|^{m-4}\left<D^2v\nabla v,\nabla v\right>,\nonumber
	\end{align}
	we obtain
	\begin{align}\label{2-2}
		\Delta v=z^{1-\frac{m}{2}}\Delta_mv-\frac{m-2}{2}\frac{\left<\nabla z,\nabla v\right>}{z},\quad\mbox{on $\{z>0\}$}.
	\end{align}
	A routine computation yields that
	\begin{align}\label{2-3}
		\left(\Delta v\right)^2&=z^{2-m}\left(\Delta_mv\right)^2
		-(m-2)z^{-\frac{m}{2}}\Delta_mv\left<\nabla z,\nabla v\right>
		\nonumber\\&\quad
		+\frac{(m-2)^2}{4}\frac{\left<\nabla z,\nabla v\right>^2}{z^2},
	\end{align}
	\begin{align}
		\nabla\Delta v&=
		z^{1-\frac{m}{2}}\nabla\Delta_mv-\frac{m-2}{2}
		z^{-\frac{m}{2}}\Delta_mv\nabla z\nonumber\\&\quad
		+\frac{m-2}{2}\frac{\left<\nabla z,\nabla v\right>\nabla z}{z^2}
		-\frac{m-2}{2}\frac{\nabla\left<\nabla z,\nabla v\right>}{z},\nonumber
	\end{align}
	and
	\begin{align}\label{2-4-4-4}
		\left<\nabla\Delta v,\nabla v\right>&=
		z^{1-\frac{m}{2}}\left<\nabla\Delta_mv,\nabla v\right>-\frac{m-2}{2}
		z^{-\frac{m}{2}}\Delta_mv\left<\nabla z,\nabla v\right>\nonumber\\&\quad
		+\frac{m-2}{2}\frac{\left<\nabla z,\nabla v\right>^2}{z^2}
		-\frac{m-2}{2}\frac{\left<\nabla\left<\nabla z,\nabla v\right>,\nabla v\right>}{z}.
	\end{align}
	Noting that $$\nabla\left<\nabla z,\nabla v\right>=D^2z\nabla v+D^2v\nabla z,$$
	and
	$$\left<D^2v\nabla z,\nabla v\right>=\left<D^2v\nabla v,\nabla z\right>=\frac{1}{2}|\nabla z|^2,$$
	we get
	\begin{align}\label{2-4-4}
		\left<\nabla\left<\nabla z,\nabla v\right>,\nabla v\right>=
		\left<D^2z\nabla v,\nabla v\right>+\frac{1}{2}|\nabla z|^2.
	\end{align}
	Combining (\ref{2-4-4}) with (\ref{2-4-4-4}), we have
	\begin{align}\label{2-4}
		\left<\nabla\Delta v,\nabla v\right>&=
		z^{1-\frac{m}{2}}\left<\nabla\Delta_mv,\nabla v\right>
		-\frac{m-2}{2}z^{-\frac{m}{2}}\Delta_mv\left<\nabla z,\nabla v\right>
		\nonumber\\&\quad
		+\frac{m-2}{2}\frac{\left<\nabla z,\nabla v\right>^2}{z^2}
		-\frac{m-2}{2}\frac{\left<D^2z\nabla v,\nabla v\right>}{z}
		\nonumber\\&\quad
		-\frac{m-2}{4}\frac{|\nabla z|^2}{z}.
	\end{align}
	By the B$\ddot{\rm o}$chner formula, we have
	\begin{align}\label{2-5}
		\frac{1}{2}\Delta |\nabla v|^2&=|D^2v|^2+\left<\nabla \Delta v, \nabla v\right>\nonumber\\
		&\geq\frac{1}{N}(\Delta v)^2+\left<\nabla \Delta v, \nabla v\right>.
	\end{align}
	Replacing \eqref{2-3} and \eqref{2-4} into (\ref{2-5}), we deduce
	\begin{align}
		\frac{1}{2}\Delta z&\geq-\frac{m-2}{2}\frac{\left<D^2z\nabla v,\nabla v\right>}{z}
		+\frac{1}{N}z^{2-m}\left(\Delta_mv\right)^2
		\nonumber\\&\quad
		+z^{1-\frac{m}{2}}\left<\nabla\Delta_mv,\nabla v\right>
		-\frac{(N+2)(m-2)}{2N}z^{-\frac{m}{2}}\Delta_mv\left<\nabla z,\nabla v\right>
		\nonumber\\&\quad
		+\frac{(2N+m-2)(m-2)}{4N}\frac{\left<\nabla z,\nabla v\right>^2}{z^2}
		-\frac{m-2}{4}\frac{|\nabla z|^2}{z}.\nonumber
	\end{align}
	The above inequality can be rewritten as
	\begin{align}
		&\quad\frac{1}{2}\mathscr A_v(z)+\frac{1}{N}z^{2-m}\left(\Delta_mv\right)^2
		+z^{1-\frac{m}{2}}\left<\nabla\Delta_mv,\nabla v\right>
		\nonumber\\&
		\leq \frac{(N+2)(m-2)}{2N}z^{-\frac{m}{2}}\Delta_mv\left<\nabla z,\nabla v\right>+\frac{m-2}{4}\frac{|\nabla z|^2}{z}
		\nonumber\\&\quad
		-\frac{(2N+m-2)(m-2)}{4N}\frac{\left<\nabla z,\nabla v\right>^2}{z^2},\nonumber
	\end{align}
which yields (\ref{L-1}).
\end{proof}

The following Bernstein estimate for solutions of (\ref{1-1}) is essential in the proofs of Theorems \ref{thm1}, \ref{thm2} and \ref{thm3}.
\begin{lemma}\label{lem2}
\rm{ Assume that $u$ is a $C^1$ solution of (\ref{1-1}) in a domain $\Omega$, with $m>1$ and $M,p,q$ arbitrary real numbers. Let $z=|\nabla u|^2$. Then for any $0<a\leq\frac{1}{N}$ and $0<b\leq\frac{M^2}{N}$, there exists a positive constant $c_1=c_1(N,M,m,q,a,b)$ such that
\begin{align}\label{2-1}
&\frac{1}{2}\mathscr A_u(z)+au^{2p}z^{2-m}
+\frac{2M}{N}|u|^{p-1}uz^{\frac{q}{2}-m+2}
\nonumber\\&\qquad
+bz^{q-m+2}-p|u|^{p-1}z^{2-\frac{m}{2}}\leq c_1\frac{\left|\nabla z\right|^2}{z},\quad\mbox{on $\{z>0\}$}.
\end{align}
}
\end{lemma}
\rm \begin{proof}
By (\ref{1-1}), we have
\begin{align}
z^{2-m}\left(\Delta_mu\right)^2=u^{2p}z^{2-m}+2M|u|^{p-1}uz^{\frac{q}{2}-m+2}+M^2z^{q-m+2},\nonumber
\end{align}
\begin{align}
z^{1-\frac{m}{2}}\left<\nabla\Delta_mu,\nabla u\right>
&=-p|u|^{p-1}z^{2-\frac{m}{2}}-\frac{Mq}{2}z^{\frac{q-m}{2}}\left<\nabla z,\nabla u\right>,\nonumber
\end{align}
and
\begin{align}
z^{-\frac{m}{2}}\Delta_mu\left<\nabla z,\nabla u\right>=-|u|^{p-1}uz^{-\frac{m}{2}}\left<\nabla z,\nabla u\right>-Mz^{\frac{q-m}{2}}\left<\nabla z,\nabla u\right>.\nonumber
\end{align}
Inserting these identities into (\ref{L-1}), we arrive
\begin{align}\label{2-9}
&\quad\frac{1}{2}\mathscr A_u(z)+\frac{1}{N}u^{2p}z^{2-m}+\frac{2M}{N}|u|^{p-1}uz^{\frac{q}{2}-m+2}
\nonumber\\&\quad
+\frac{M^2}{N}z^{q-m+2}-p|u|^{p-1}z^{2-\frac{m}{2}}
\nonumber\\ &\leq
-\frac{(N+2)(m-2)}{2N}|u|^{p-1}uz^{-\frac{m}{2}}\left<\nabla z,\nabla u\right>
\nonumber\\&\quad
+\left(\frac{Mq}{2}-\frac{M(N+2)(m-2)}{2N}\right)z^{\frac{q-m}{2}}\left<\nabla z,\nabla u\right>+\frac{m-2}{4}\frac{|\nabla z|^2}{z}
\nonumber\\&\quad
-\frac{(2N+m-2)(m-2)}{4N}\frac{\left<\nabla z,\nabla u\right>^2}{z^2}
,\quad\mbox{on $\{z>0\}$}.
\end{align}
Next we estimate each term in the right-hand side of (\ref{2-9}).
By Cauchy-Schwartz inequality and then, thanks to Young inequality,  we have
for any $\varepsilon,\varepsilon'>0$
$$|u|^{p-1}uz^{-\frac{m}{2}}\left|\left<\nabla z,\nabla u\right>\right|\leq \varepsilon u^{2p}z^{2-m}+\frac{1}{ 4\varepsilon}\frac{|\nabla z|^2}{z},$$
and
$$z^{\frac{q-m}{2}}\left|\left<\nabla z,\nabla u\right>\right|\leq \varepsilon ' z^{q-m+2}+\frac{1}{4\varepsilon'}\frac{|\nabla z|^2}{z}.$$
Note also that
$$\frac{\left<\nabla z,\nabla u\right>^2}{z^2}\leq \frac{|\nabla z|^2}{z}.$$
Let $\varepsilon_1:=\frac{(N+2)|m-2|}{2N}\varepsilon$ and $\varepsilon_2:=\left|\frac{Mq}{2}-\frac{M(N+2)(m-2)}{2N}\right|\varepsilon'$.
We infer that
\begin{align}
&\frac{1}{2}\mathscr A_u(z)+\left(\frac{1}{N}-\varepsilon_1\right)u^{2p}z^{2-m}
+\frac{2M}{N}|u|^{p-1}uz^{\frac{q}{2}-m+2}
\nonumber\\&\quad
+\left(\frac{M^2}{N}-\varepsilon_2\right)
z^{q-m+2}-p|u|^{p-1}z^{2-\frac{m}{2}}\leq c_1\frac{\left|\nabla z\right|^2}{z},
\nonumber
\end{align}
where $c_1=c_1(N,m,\varepsilon_1,\varepsilon_2)>0$. Set $a=\frac{1}{N}-\varepsilon_1$ and $b=\frac{M^2}{N}-\varepsilon_2$. Taking $\varepsilon_1$ and $\varepsilon_2$ small enough such that $a,b>0$, then (\ref{2-1}) follows.
\end{proof}

Now we step into the proof of Theorem \ref{thm1}.
\rm \begin{proof}[\rm\textbf{Proof of Theorem \ref{thm1}}]
Let $u$ be a positive solution of \eqref{1-1}. Consider the following change of variables
\begin{eqnarray}\label{def-v}
u(x)=\alpha^{\frac{m}{p-m+1}}v(y), \quad   y=\alpha x, \quad x\in\Omega,
\end{eqnarray}
with $\alpha=M^{-\frac{p-m+1}{(p+1)q-mp}}.$

Then  $|\nabla v|=|\nabla_yv|=\alpha^{-\frac{p+1}{p-m+1}}|\nabla u|$ and $\Delta_m v=\alpha^{-\frac{mp}{p-m+1}}\Delta_m u$ so that
$v$ is a positive $C^1$ solution of
\begin{equation}\label{2}
-\Delta_mv=|v|^{p-1}v+\left|\nabla v\right|^q  \quad\text{in } \Omega_\alpha,
\end{equation}
 where $ \Omega_\alpha:=\{y\in\mathbb R^N: \, y=\alpha x, \,\, x\in\Omega\}$.

Let  $z=|\nabla v|^2$, so that  \eqref{2-1} becomes
\begin{align*}
&\frac{1}{2}\mathscr A_v(z)+av^{2p}z^{2-m}
+\frac{2}{N}|v|^{p-1}vz^{\frac{q}{2}-m+2}
\nonumber\\&\quad
+bz^{q-m+2}-p|v|^{p-1}z^{2-\frac{m}{2}}\leq c_1\frac{\left|\nabla z\right|^2}{z},\quad\mbox{on $\{z>0\}$},
\end{align*}
indeed $v$ is a positive solution of \eqref{1-1} with $M=1$. In turn
\begin{align}\label{2-1-1}
&\frac{1}{2}\mathscr A_v(z)+av^{2p}z^{2-m}
+bz^{q-m+2}\nonumber\\&\quad
-p|v|^{p-1}z^{2-\frac{m}{2}}\leq c_1\frac{\left|\nabla z\right|^2}{z},\quad\mbox{on $\{z>0\}$},
\end{align}
with $0<a,b\le\frac{1}{N}$ as in  \eqref{2-1} and $c_1=c_1(N,m,q,a,b)$.

Suppose $q>\frac{mp}{p+1}$.  In this case, it immediately follows that $q-m+2>1$ by conditions assumed on $p$.
 By Young inequality with exponents $2p/(p-1)$ and $2p/(p+1)$, for  $\varepsilon_3>0$, we have
$$p|v|^{p-1}z^{2-\frac{m}{2}}=p|v|^{p-1}z^{\frac{(2-m)(p-1)}{2p}}z^{1+\frac{2-m}{2p}}
\leq\varepsilon_3v^{2p}z^{2-m}+c_2z^{\frac{2p+2-m}{p+1}}.$$
Since $2p>m-2$ and $q(p+1)>mp$, a further application of Young inequality with exponents $(q-m+2)(p+1)/(2p+2-m)$ and its conjugate, gives,  for $\varepsilon_4>0$,
$$c_2z^{\frac{2p+2-m}{p+1}}\leq\varepsilon_4 z^{q-m+2}+c_3,$$
where $c_2=c_2(p,\varepsilon_3)>0$ and $c_3=c_3(m,p,q,c_2,\varepsilon_4)>0$.
Hence by (\ref{2-1-1}),
\begin{align}
&\frac{1}{2}\mathscr A_v(z)+A_1v^{2p}z^{2-m}
+A_2z^{q-m+2}
\leq c_1\frac{\left|\nabla z\right|^2}{z}+c_3,\nonumber
\end{align}
where $A_1=a-\varepsilon_3$ and $A_2=b-\varepsilon_4$.
Taking $\varepsilon_3$ and $\varepsilon_4$ small enough such that $A_1,A_2>0$, then
\begin{align}
&\frac{1}{2}\mathscr A_v(z)+A_2z^{q-m+2}
\leq c_1\frac{\left|\nabla z\right|^2}{z}+c_3.\nonumber
\end{align}
Letting $\widetilde{z}=\left(z-\left(\frac{c_3}{A_2}\right)^{\frac{1}{q-m+2}}\right)_+$,  thus $z\ge \tilde z$,
	and being $q-m+2>1$, we obtain
\begin{align}
&\frac{1}{2}\mathscr A_v(\widetilde{z})+A_2\widetilde{z}^{q-m+2}
\leq c_1\frac{\left|\nabla \widetilde{z}\right|^2}{\widetilde{z}},
\quad\mbox{on $\left\{z>(\frac{c_3}{A_2})^{\frac{1}{q-m+2}}\right\}$
}.\nonumber
\end{align}

Using Lemma \ref{lem1}, we derive
$$\widetilde{z}(y)\leq c_4\left({\rm dist}\left(y,\partial\Omega_\alpha\right)\right)^{-
\frac{2}{q-m+1}},$$
where $c_4=c_4(m,q,c_1,A_2)>0$, and being $\tilde z= |\nabla v(y)|^2- c$, $c>0$, using that $(a+b)^{1/2}\le a^{1/2}+b^{1/2}$,  it follows that
\begin{align}\label{3-1}
|\nabla v(y)|\leq c_4'\left(1+\left({\rm dist}\left(y,\partial\Omega_\alpha\right)\right)^{-
\frac{1}{q-m+1}}\right),\quad y\in\Omega_\alpha.
\end{align}
In view of the change of variables \eqref{def-v}, we finally obtain
 \eqref{1-2}.

Now consider the case $\Omega=\mathbb R^N$ and assume that $u$ is a positive solution  of \eqref{1-1} in $\mathbb R^N$. Fix $y\in\mathbb R^N$ such that $|y|<2n$. Using (\ref{3-1}) with $\Omega_\alpha=B_{2n}(0)$, we see
$$|\nabla v(y)|\leq c_4'\left(1+\left(2n-|y|\right)^{-\frac{1}{q-m+1}}\right),\quad  y\in B_{2n}(0).$$
Taking $n\rightarrow\infty$ yields
\begin{align}
|\nabla v(y)|\leq c_4',\quad y\in\mathbb R^N,\nonumber
\end{align}
so that \eqref{1-3} follows immediately thanks to the change of variables.
\end{proof}

\rm \begin{proof}[\rm\textbf{Proof of Theorem \ref{thm2}}]  Let $u$ be a positive solution of \eqref{1-1}  and let $v$ be the function defined in \eqref{def-v} where now
$\Omega=\Omega_\alpha=\mathbb R^N$. If $z=|\nabla v|^2$,
	since we have $\max\left\{m-1,\frac{m}{2}\right\}<q<\frac{mp}{p+1}$, then for any $\varepsilon_5>0$, we have
$$pv^{p-1}z^{2-\frac{m}{2}}=pv^{p-1}z^{\frac{(2-m)(2q-m)}{2q}}z^{\frac{m(q-m+2)}{2q}}
\leq\varepsilon_5z^{q-m+2}+c_5v^{\frac{2q(p-1)}{2q-m}}z^{2-m},$$
where $c_5=c_5(m,p,q,\varepsilon_5)>0$.
Inserting this inequality into (\ref{2-1-1}), we obtain
\begin{align}
&\frac{1}{2}\mathscr A_v(z)+v^{2p}z^{2-m}\left(a-c_5v^{\frac{2mp-2q(p+1)}{2q-m}}\right)
+A_3z^{q-m+2}\leq c_1\frac{\left|\nabla z\right|^2}{z},\nonumber
\end{align}
where $A_3=b-\varepsilon_5$ with  $\varepsilon_5$ small enough such that $A_3>0$.
If $$\max v\leq c_{N,m,p,q}:=\left(\frac{a}{c_5}\right)^{\frac{2q-m}{2mp-2q(p+1)}},$$
which is equivalent to (\ref{1-4}) by virtue of \eqref{def-v}, we get
\begin{align}
&\frac{1}{2}\mathscr A_v(z)
+A_3z^{q-m+2}\leq c_1\frac{\left|\nabla z\right|^2}{z}.\nonumber
\end{align}
From Lemma \ref{lem1},  applied with $\xi=q-m+2>1$,  we conclude that $z\equiv0$ in $\mathbb R^N$, in turn $v$ is
identically constant and thus $v\equiv0$ in $\mathbb R^N$ from the equation \eqref{2}.
\end{proof}

\rm \begin{proof}[\rm\textbf{Proof of Theorem \ref{thm3}}]Let $u$ be a positive solution of \eqref{1-1} in $\Omega$ and let $q=\frac{mp}{p+1}$ by assumption.
Consider the auxiliary function $\Phi$
defined  for $Z>0$ by
$$\Phi(Z)=u^pZ^{2-m}+BZ^{q-m+2}-\sqrt{\frac{p}{a}}u^{\frac{p-1}{2}}Z^{2-\frac{m}{2}}.$$
In particular $\Phi(Z)=Z^{2-m}\psi (Z)$ where
$$\psi (Z)=u^p+BZ^{\frac{mp}{p+1}}-\sqrt{\frac{p}{a}}u^{\frac{p-1}{2}}Z^{\frac{m}{2}}$$
with
$$\psi (0)=u^p>0\quad\text{and  } \psi'(Z)=\frac{mBp}{p+1}Z^{\frac m2 -1}\biggl[Z^{\frac{m(p-1)}{2(p+1)}}-\frac{p+1}{2\sqrt{ap}}u^{\frac{p-1}{2}}\biggr],$$
so that $\psi(Z)$ achieves its minimum at
$$Z_0=\left(\frac{p+1}{2B\sqrt{ap}}\right)^{\frac{2(p+1)}{m(p-1)}}u^{\frac{p+1}{m}}>0,$$
 and
$$\psi(Z)\ge \psi(Z_0)=\left[1-\frac{p-1}{(4ap)^{\frac{p}{p-1}}}\left(\frac{p+1}{B}\right)^{\frac{p+1}{p-1}}\right]u^p.$$
Denoting
$$M_+=\frac{(p+1)(p-1)^{\frac{p-1}{p+1}}}{(4ap)^{\frac{p}{p+1}}}>0,$$
we obtain that if $B\geq M_+$, then {$\psi(Z_0)\geq0$ yielding } $\psi(Z)\geq0$ for all $Z>0$ and consequently $\Phi(Z)\geq0$ for all $Z>0$.

Now consider inequality (\ref{2-1}) for $u$ positive in the set where $|\nabla u|\neq0$.  If $M>aNM_+$, we have
\begin{align}
&\frac{1}{2}\mathscr A_u(z)+a\left(u^{p}z^{1-\frac{m}{2}}+M_+z^{1+\frac{q-m}{2}}\right)^2
\nonumber\\&\qquad
+(b-aM_+^2)z^{q-m+2}-pu^{p-1}z^{2-\frac{m}{2}}\leq c_1\frac{\left|\nabla z\right|^2}{z}.\nonumber
\end{align}
We claim that
\begin{equation}\label{ineq_use}a\left(u^{p}z^{1-\frac{m}{2}}+M_+z^{1+\frac{q-m}{2}}\right)^2-pu^{p-1}z^{2-\frac{m}{2}}\ge0.\end{equation}
 Indeed,   noting  that for any $B\geq M_+$, we have
\begin{align}
&\quad\left(u^p|\nabla u|^{2-m}+B|\nabla u|^{q-m+2}\right)^2-\frac{p}{a}u^{p-1}|\nabla u|^{4-m}
\nonumber\\&
=\left(u^p|\nabla u|^{2-m}+B|\nabla u|^{q-m+2}+\sqrt{\frac{p}{a}}u^{\frac{p-1}{2}}|\nabla u|^{2-\frac{m}{2}}\right) \cdot\Phi(|\nabla u|)\ge0,
\nonumber
\end{align}
where $0<a\le\frac{1}{N}$, then   \eqref{ineq_use}  immediately follows choosing $B=M_+$,  being  $z=|\nabla u|^2$.

Consequently,
\begin{align}
&\frac{1}{2}\mathscr A_u(z)+(b-aM_+^2)
z^{q-m+2}\leq c_1\frac{\left|\nabla z\right|^2}{z}.\nonumber
\end{align}
Letting $aM_+^2<b\leq\frac{M^2}{N}$ and using again Lemma \ref{lem1},  with $\xi=q-m+2>1$, we obtain
\begin{equation}
\left|\nabla u(x)\right|\leq c_{N,M,a,m,p,q}
\left({\rm dist}\left(x,\partial\Omega\right)\right)^{-
\frac{1}{q-m+1}},\nonumber
\end{equation}
which is exactly (\ref{1-5-5}) via $q=\frac{mp}{p+1}$.
\end{proof}

\section{Proof of Theorems \ref{thm4}, \ref{thm5} and \ref{thm0}}

\begin{proposition}\label{prop}
\rm{
\textbf{(A)}. Let $M\geq0$, $m>1$, $q\geq0$ and either $N\leq m$ and $p>0$ or $N>m$ and $0<p\leq\frac{N(m-1)}{N-m}$. Then, there exist no positive solutions of (\ref{1-1}) in $\mathbb R^N\setminus \overline{B}_R$ for $R>0$.

\noindent \textbf{(B).} Let  $M>0$, $m>1$, $N>1$, $p\geq0$ and $m-1<q\leq\frac{N(m-1)}{N-1}$, then there exist no positive radial solutions of (\ref{1-1}) in $\mathbb R^N\setminus \overline{B}_R$ for $R\ge0$.

\noindent \textbf{(C).} Let  $N>m$, $m>1$, $M\geq0$, $q\geq 0$, $p>\frac{N(m-1)}{N-m}$ and let $u=u(|x|)=u(r)$ be a positive radial solution of (\ref{1-1}) in $\mathbb R^N\setminus \overline{B}_R$. Then, there exists $\rho> R$ such that
\begin{align}\label{2--1}
 u(r)\leq c_0r^{-\frac{m}{p-m+1}},\quad  r>\rho,
\end{align}
with $c_0=\left[2N\left(1-2^{-\frac{m}{m-1}}\right)^{-(m-1)}\left(\frac{m}{p-m+1}\right)^{m-1}\right]^{\frac{1}{p-m+1}}$ and
\begin{align}\label{2--2}
\left|u_r(r)\right|\leq c_0\frac{N-m}{m-1}r^{-\frac{p+1}{p-m+1}},\quad  r>\rho.
\end{align}

\noindent \textbf{(D).} Let $N>1$, $m>1$, $M>0$, $p\geq0$, $q>\frac{N(m-1)}{N-1}$,
and let $u(x)=u(r)$ be a positive radial solution of (\ref{1-1}) in $\mathbb R^N\setminus \overline{B}_R$. There exists $\rho>2R$
such that
\begin{align}\label{2--3}
\left|u_r(r)\right|\leq c_1r^{-\frac{1}{q-m+1}},\quad r>\frac{\rho}{2},
\end{align}
with $c_1=\left(\frac{q(N-1)-N(m-1)}{M(q-m+1)}\right)^{\frac{1}{q-m+1}}$. Moreover, if $\frac{N(m-1)}{N-1}<q<m$
\begin{align}\label{2--4}
u(r)\leq c_1\frac{q-m+1}{m-q}r^
{-\frac{m-q}{q-m+1}},\quad r>\frac{\rho}{2}.
\end{align}
}\end{proposition}

\rm \begin{proof}\textbf{
(A):} When $M\geq0$,  every solution  $u$ of  (\ref{1-1}) satisfies the inequality
	$$
		-\Delta_mu\geq|u|^{p-1}u,\quad{\rm in}\;\mathbb R^N\setminus \overline{B}_R.$$
Then, assertion (A) follows by Theorems 3.3 (iii) and 3.4 (ii) of  \cite{Pohozaev} and  Theorem I$'$ in \cite{SZ}.
	
\textbf{(B):} Let $u$ be a radial positive solution  of \eqref{1-1} in  $\mathbb R^N\setminus \overline{B}_R$, $R\ge0$. Thus, $u=u(r)=u(|x|)$ satisfies \eqref{1-1} in the radial form,
that is
	\begin{align}\label{p-1}
		-r^{1-N}\left(r^{N-1}|u_r|^{m-2}u_r\right)_r=u^p+M\left| u_r\right|^q, \quad r> R.
	\end{align}	
	It follows that $r\mapsto w(r):=-r^{N-1}|u_r|^{m-2}u_r$ is strictly increasing on $(R,\infty)$, thus it admits a limit $l\in(-\infty,\infty]$. If $l\leq0$, then $u_r(r)>0$ on $(R,\infty)$. Hence $u(r)\geq{u}(s_0):=c>0$ for some $s_0>R$  and for all $r\geq s_0$, so that
		$$\left(r^{N-1}u_r^{m-1}\right)_r\leq-c^pr^{N-1}, \quad r\geq s_0,$$
 in turn, by integration form $s$ to $r$, with $s_0<s<r$, we arrive to
	$$\left(u_r(r)\right)^{m-1}\leq\frac{s^{N-1}}{r^{N-1}}\left( u_r(s)\right)^{m-1}-\frac{c^p}{N}\left(r-\frac{s^N}{r^{N-1}}\right),$$
	which implies $u_r(r)\to-\infty$, thus $u(r)\to-\infty$ as $r\to\infty$, a contradiction. Therefore, $w(r)\to l\in(0,\infty]$ as $r\to\infty$ and there exists $r_l>R$ such that $u_r(r)<0$ on $(r_l,\infty)$,
 so that $w=r^{N-1}|u_r|^{m-1}>0$ on $(r_l,\infty)$.
	By (\ref{p-1}), we have for $M>0$
	$$w_r\geq Mr^{-\frac{(N-1)(q-m+1)}{m-1}}w^{\frac{q}{m-1}},$$	
yielding
\begin{align}\label{p--2}
	\left(w^{-\frac{q-m+1}{m-1}}\right)_r\leq-\frac{q-m+1}{m-1}Mr^{-\frac{(N-1)(q-m+1)}{m-1}}.
\end{align}	
Integrating (\ref{p--2}) on $(s,r)$ with $s>r_l$, if $q=\frac{N(m-1)}{N-1}$, we obtain
	\begin{align}\label{p-2}
		w^{-\frac{1}{N-1}}(r)-w^{-\frac{1}{N-1}}(s)\leq-\frac{M}{N-1}\ln\frac{r}{s},
	\end{align}	
while if $q<\frac{N(m-1)}{N-1}$, we have
\begin{align}\label{p-3}
	&w^{-\frac{q-m+1}{m-1}}(r)-w^{-\frac{q-m+1}{m-1}}(s)
	\nonumber\\ \leq&-\frac{M(q-m+1)}{N(m-1)-q(N-1)}
	\left(r^{\frac{N(m-1)-q(N-1)}{m-1}}-s^{\frac{N(m-1)-q(N-1)}{m-1}}\right).
\end{align}	
Letting $r\to\infty$, we obtain that both right-hand sides of (\ref{p-2}) and (\ref{p-3}) tend to $-\infty$, being $N(m-1)-q(N-1)>0$, namely
$$w^{-\frac{1}{N-1}}(r),\, w^{-\frac{q-m+1}{m-1}}(r)\to \-\infty\quad\text{as }r\to\infty.$$
This contradicts $\lim\limits_{r\to\infty}w(r)=l>0$, concluding the proof of (B).

\textbf{(C):} Let $u(x)=u(r)$ be a positive radial solution of (\ref{1-1}) in $\mathbb R^N\setminus \overline{B}_R$. Arguing as in \textbf{(B)}, there exists $r_l> R$ such
 that $u_r(r)<0$ on $(r_l,\infty)$. By (\ref{p-1}), being $M\ge0$, we have for $r>\rho_1:=2r_l$
	\begin{align}
		r^{N-1}\left| u_r(r)\right|^{m-1}\geq\int_{\frac{r}{2}}^r\tau^{N-1}{u}^p(\tau)d\tau\geq\frac{r^N{u}^p(r)}{N}\biggl (1-\frac1{2^N}\biggr)\ge \frac{r^N{u}^p(r)}{2N} ,\nonumber
	\end{align}
 yielding
\begin{align}\label{p-4}
	\left({u}^{-\frac{p-m+1}{m-1}}\right)_r\geq\frac{p-m+1}{m-1}\left(\frac{r}{2N}\right)^{\frac{1}{m-1}}.
\end{align}
Integrating (\ref{p-4}) on $\left(\frac{r}{2},r\right)$ we obtain
\begin{align}\label{p-5}
	 u(r)\leq c_0r^{-\frac{m}{p-m+1}},\quad  r>2\rho_1,
\end{align}
	with $c_0=(2N)^{\frac{1}{p-m+1}}\bigl(1-2^{-\frac{m}{m-1}}\bigr)^{-\frac{m-1}{p-m+1}}\left(\frac{m}{p-m+1}\right)^{\frac{m-1}{p-m+1}}$, which yields (\ref{2--1}).
	
To prove \eqref{2--2}, we set $v(t)=u(t^{-\frac{m-1}{N-m}})$ with  $t\in\left(0,\rho_1^{-\frac{N-m}{m-1}}\right)$. By (\ref{p-5}) we see that $v(t)\to0$ as $t\to0^+$.
By \eqref{p-1}, using that $v_t(t)=-\frac{m-1}{N-m}u_r (r)r^{\frac{N-1}{m-1}}$, and $r=t^{-\frac{m-1}{N-m}}$, we obtain
	\begin{align}
		v_{tt}(t)&=\frac{m-1}{(N-m)^2}r^{\frac{2(N-1)}{m-1}}\left[(m-1) u_{rr}+\frac{N-1}{r}u_r\right]
		\nonumber\\&
		=\frac{m-1}{(N-m)^2}r^{\frac{(3-m)(N-1)}{m-1}}| u_r|^{2-m}\left(r^{N-1}|u_r|^{m-2} u_r\right)_r\leq0.\nonumber
	\end{align}
	Using mean value theorem in $(0,t)$, we derive, being  $v_t$ is increasing since $u_r<0$,
	$$v_t(t)\leq\frac{v(t)}{t},$$
 so that, replacing the expression of $v_t$, we obtain the following 
$$\left|u_r(r)\right|\leq\frac{N-m}{m-1}t^{\frac{N-1}{N-m}}\frac{v(t)}{t}=\frac{N-m}{m-1} \frac{v(t)}{t^{-\frac{m-1}{N-m}}}=\frac{N-m}{m-1}\frac{u(r)}{r},\qquad
	 r>2\rho_1,$$
	 so that, using \eqref{2--1} with $\rho=2\rho_1$, then \eqref{2--2}  follows immediately.
	
	\textbf{(D):} Let $u$ be a radial positive solution  of \eqref{1-1} in  $\mathbb R^N\setminus \overline{B}_R$, $R>0$. Arguing as in the first part of \textbf{(B)}, but now assuming $q>\frac{N(m-1)}{N-1}$,  inequality  \eqref{p-3} is still valid, so that letting $r\to\infty$ on  both sides of (\ref{p-3}), we obtain
that there exists $\rho$ such that for all $s>\frac{\rho}{2}$
	$$l^{-\frac{q-m+1}{m-1}}-w^{-\frac{q-m+1}{m-1}}(s)\leq-\frac{M(q-m+1)}{q(N-1)-N(m-1)}s^{-\frac{q(N-1)-N(m-1)}{m-1}},$$
hence
	$$w(s)\leq\left(\frac{q(N-1)-N(m-1)}{M(q-m+1)}\right)^{\frac{m-1}{q-m+1}}s^{\frac{q(N-1)-N(m-1)}{q-m+1}},\quad s>\frac{\rho}{2},$$
 thus, form $w(r)=r^{N-1}|w_r(r)|^{m-1}$ we get
	$$\left|u_r(r)\right|\leq\left(\frac{q(N-1)-N(m-1)}{M(q-m+1)}\right)^{\frac{1}{q-m+1}}r^{-\frac{1}{q-m+1}},\quad r>\frac{\rho}{2},
	$$
	which yields (\ref{2--3}). Then (\ref{2--4}) follows by integrating (\ref{2--3}) from $r$ to $\infty$.
\end{proof}

\rm \begin{proof}[\rm\textbf{Proof of Theorem \ref{thm4}}]
Let $u$ be a positive supersolution of (\ref{1-1}) in $\overline{B}_R^c$ for  some $R>0$. 
By Proposition \ref{prop} (A), we know that when $M\geq0$, the result is valid,  even in a larger range for $p$. Thus, let us deal with  the remaining case $M<0$ and $N=m$ with $p>m-1$ or $N>m$ with $m-1<p<\frac{N(m-1)}{N-m}$.

Setting $u=v^\sigma$ with $\sigma>1$, we obtain
\begin{align}
-\Delta_mv&\geq(\sigma-1)(m-1)\frac{|\nabla v|^m}{v}+\sigma^{1-m}v^{m+\sigma(p-m+1)-1}
\nonumber\\&\quad
+M\sigma^{q-m+1}v^{(\sigma-1)(q-m+1)}
|\nabla v|^q,\nonumber
\end{align}
and then setting $z=|\nabla v|^m$ yields
\begin{align}\label{4-1}
-\Delta_mv\geq\sigma^{1-m}\frac{\Psi(z)}{v},
\end{align}
where
$$\Psi(z)=\sigma^{m-1}(\sigma-1)(m-1)z+M\sigma^{q}v^{(\sigma-1)(q-m+1)+1}z^{\frac{q}{m}}+v^{m+\sigma(p-m+1)}.$$
Since $q=\frac{mp}{p+1}$, it is easy to see that $\Psi(z)$ achieves its minimum at
$$z_0=\left(\frac{|M|p\sigma^{1-\frac{m}{p+1}}}{(\sigma-1)(m-1)(p+1)}\right)^{p+1}v^{m+\sigma(p-m+1)},$$
and
\begin{align}\label{4-2}
\Psi(z_0)=\left[1-\left(\frac{|M|}{p+1}\right)^{p+1}\left(\frac{\sigma p}{(\sigma-1)(m-1)}\right)^p\right]v^{m+\sigma(p-m+1)}.
\end{align}
For the case of $N>m$, we choose $\sigma$ such that
$$m+\sigma(p-m+1)-1=\frac{N(m-1)}{N-m},$$
namely
$$\sigma=\frac{m(m-1)}{(N-m)(p-m+1)},$$
in turn  $\sigma>1$ by $p<\frac{N(m-1)}{N-m}$ and
\begin{align}
 \Psi(z)\ge\Psi(z_0)=\left[1-\left(\frac{|M|}{p+1}\right)^{p+1}\left(
\frac{mp}{N(m-1)-p(N-m)}\right)^p\right]v^{\frac{N(m-1)}{N-m}+1}.
\nonumber
\end{align}
We derive that if $|M|<\mu^\ast(N)$, where $\mu^\ast(N)$ is given in \eqref{1-6},
then  inequality \eqref {4-1} gives
\begin{equation}\label{no_ext}-\Delta_mv\geq\delta v^{\frac{N(m-1)}{N-m}} \quad\text{in } \mathbb R^N\setminus \overline{B}_R,\end{equation}
    for some $\delta>0$.  Hence,  Proposition \ref{prop} (A)  yields the required contradiction,   since no positive solutions of \eqref{no_ext} can exist in  exterior domains of $\mathbb R^N$.

If $N=m$, for a fixed $\sigma>1$,  if
\begin{align}\label{4-4}
	|M|<(p+1)\left(\frac{(\sigma-1)(m-1)}{\sigma p}\right)^{\frac{p}{p+1}}:=\mu^*_m,
\end{align}
then,  from \eqref{4-1} and \eqref{4-2}, we have 
$$-\Delta_mv\geq\delta v^{m+\sigma(p-m+1)-1} \quad\text{in } \mathbb R^N\setminus \overline{B}_R$$
for some $\delta>0$. Since $m+\sigma(p-m+1)-1>0$, then the result follows immediately from Proposition \ref{prop} (A).   In particular,  $$\mu^*_m \to \mu^\ast(m)=(p+1)\biggr(\frac{m-1}p\biggr)^{\frac{p}{p+1}}\quad\text{as } \sigma\to\infty,$$
thus, choosing $\sigma$ large  enough, condition $M>-\mu^*(N)$ holds also for $N=m$.
\end{proof}

\rm \begin{proof}[\rm\textbf{Proof of Theorem \ref{thm5}}] We perform the proof by contradiction argument. Let us assume that there exists a positive supersolution $u$ of \eqref{1-1} satisfying (\ref{1-7}).
Without loss of generality,  let us assume that $u>1$ in $\Omega$, otherwise, we could replace $\Omega$ with the set
$\{u>1\}$.
Take   $v=\log u$, so that $v$ is positive being $u>1$.
 By $q=\frac{mp}{p+1}$, we obtain
\begin{align}\label{4-7}
-\Delta_mv\geq F\left(|\nabla v|^m\right),
\end{align}
where
$$F(X)=(m-1)X+e^{(p-m+1)v}+Me^{(q-m+1)v}X^{\frac{p}{p+1}}.$$
Obviously  $F(X)>0$ for any $X\geq0$ when $M\geq0$. On the other hand, in  the case  $M<0$,
it is not hard to see that $F(X)$ achieves its minimum at
$$X_0=\left(\frac{|M|p}{(m-1)(p+1)}\right)^{p+1}e^{(p-m+1)v},$$
and
$$
F(X)\ge F(X_0)=\left[1-\left(\frac{p}{m-1}\right)^p\left(\frac{|M|}{p+1}\right)^{p+1}\right]e^{(p-m+1)v}
$$
for all $X\ge0$. Therefore, if
\begin{align}
|M|\leq(p+1)\left(\frac{m-1}{p}\right)^{\frac{p}{p+1}}=\mu^\ast(m),
\end{align}
where $\mu^\ast$ is as in \eqref{1-6}, then $F(X_0)\ge0$, so that
we see that $v$ solves
\begin{align}\label{4-7-7}
 \begin{cases}
 -\Delta_mv\geq0,\quad{\rm in}\;\,\Omega,\\
 \lim\limits_{{\rm dist}\left(x,\partial\Omega\right)\rightarrow0}v(x)=\infty.
 \end{cases}
\end{align}
Clearly, when $\Omega$ is bounded, $v$ is  larger than  the $m$-harmonic  function with any boundary value $k>0$. Letting $k\rightarrow\infty$ we derive a contradiction.

When $\Omega$ is an exterior domain, namely $\Omega=\mathbb R^N\setminus \overline{B}_R$, so that  $\Omega^c=B_R$, we may assume $B_{R_1}\subset\Omega^c\subset B_{R_2}$ for some $R_2>R_1>0$.
Define
	$$d=\frac{(N-1)(R_2-R_1)}{(m-1)R_1}+1>0$$
and
	$$w(x)=\left(R_2-|x|\right)^d,\quad{\rm in}\;\, B_{R_2}\backslash{\Omega^c}.$$
It holds
	$$-\Delta_mw=d^{m-1}\left(R_2-|x|\right)^{d(m-1)-m}\biggl[ (N-1)\frac{R_2-|x|}{|x|}-(d-1)(m-1)\biggr],$$
	thus, the choice of $d$ and the decreasing monotonicity of $(R_2-y)/y$ in $(R_1,R_2)$, gives
	 that $w$ is a solution of
	\begin{align}\label{4-7-8}
		\begin{cases}
			-\Delta_mw\leq0,&{\rm in}\;\,B_{R_2}\backslash{\Omega^c},\\
			w\leq\left(R_2-R_1\right)^d,&{\rm on}\;\,\partial\Omega,\\
			w=0,&{\rm on}\;\,\partial B_{R_2}.
		\end{cases}
	\end{align}
Hence by the weak comparison principle in \cite[Lemma 2.2]{SZ}, we get $v\geq kw$ in $B_{R_2}\backslash{\Omega^c}$ for any $k>0$. Letting $k\rightarrow\infty$ we derive a contradiction once again.
\end{proof}

\rm \begin{proof}[\rm\textbf{Proof of Theorem \ref{thm0}}]
Let $u\in C^2(\Omega\setminus\{0\})$ be a positive solution of  \eqref{1-1} in $\Omega\setminus\{0\}$. Let $\bar{B_1}\subset\Omega$. By \cite[Theorem 1.1]{Bidaut},
$$
u^{m-1}\in \mathcal M^{\frac{N}{N-m}}(B_1),\quad |\nabla u|^{m-1}\in \mathcal M^{\frac{N}{N-1}}(B_1),
$$
where $\mathcal M^r=L^{r,\infty}$ denotes the Marcinkiewicz space or Lorentz space of index $(r,\infty)$.
In order to fit with Serrin's formalism, we write (\ref{1-1}) as
$$
-\Delta_mu=Du^{m-1}+E\left|\nabla u\right|^{m-1},
$$
where $D=u^{p-m+1}$ and $E=M\left|\nabla u\right|^{q-m+1}$. Then
$$
D\in \mathcal M^{\frac{N(m-1)}{(N-m)(p-m+1)}}(B_1),\quad E\in \mathcal M^{\frac{N(m-1)}{(N-1)(q-m+1)}}(B_1).
$$
Since $m-1<p<\frac{N(m-1)}{N-m}$ and $m-1<q<\frac{N(m-1)}{N-1}$, we have
\begin{align}\label{00-0}
\frac{N(m-1)}{(N-m)(p-m+1)}>\frac{N}{m},\quad \frac{N(m-1)}{(N-1)(q-m+1)}>N.
\end{align}
Since $\mathcal M^r(B_1)\hookrightarrow L^{r-\delta}(B_1)$ for any $r>\delta>0$, we infer that
$$
D\in  L^{\frac{N}{m}+\delta}(B_1),\quad E\in L^{N+\delta}(B_1).
$$
Thus $u$ verifies the Harnack inequality in $B_1\backslash\{0\}$ by \cite[Theorem 5]{Serrin}. This implies that
\begin{align}\label{00-1}
\max_{|x|=r}u(x)\leq K\min_{|x|=r}u(x),\quad\forall r\in\left(0,1/2\right],
\end{align}
where $K>0$ depending on the norms of $D$ and $E$.

Moreover, since $u(x)=u(r)$ on $\{x:|x|=r\}$ is $m$-superharmonic when $M\geq0$, i.e.,
$
-\left(r^{N-1}| u_r|^{m-2} u_r\right)_r\geq0, 
$
 there exists some $k>0$ such that
\begin{equation}\label{k}
	u(r)\leq kr^{\frac{m-N}{m-1}}.
\end{equation}
 Indeed, by  monotonicity decreasing of $r^{N-1}| u_r|^{m-2} u_r$,   there exists $k_0>0$ such that $$r^{N-1}| u_r|^{m-2} u_r\geq-k_0,$$
	which yields
	\begin{align}\label{k-1}
		u_r\geq-k_0^{\frac{1}{m-1}}r^{\frac{1-N}{m-1}},\quad\mbox{for $r\in(0,1]$}.
	\end{align}
	Integrating (\ref{k-1}) on $(r,1)$, we obtain
	$$u(1)-u(r)\geq k(1-r^{\frac{m-N}{m-1}}),$$
	where $k=k_0^{\frac{1}{m-1}}\frac{m-1}{N-m}$. It follows that
	\begin{align}\label{k-2}
		u(r)\leq u(1)-k+kr^{\frac{m-N}{m-1}}\leq k^{\prime}r^{\frac{m-N}{m-1}},
	\end{align}
in a suitable right neighborhood of $0$,  being $m<N$,  so that \eqref{k} holds.

Combining with (\ref{00-1}), we arrive
\begin{equation*}\label{00-2}
u(x)\leq Kk|x|^{\frac{m-N}{m-1}}.
\end{equation*}
According to (\ref{00-0}) and (\ref{00-2}), we see that the function $g:=|u|^{p-1}u+M\left|\nabla u\right|^q$ satisfies the $(\phi,m)$-scaling-growth property defined by \cite[Definition 3.1]{L}, thus
the estimate on the gradient is standard and follows \cite[Lemma 3.3.2]{L}.

\end{proof}
\section{Proof of Theorem \ref{thm6}}
\rm \begin{proof}[\rm\textbf{Proof of Theorem \ref{thm6}}] Let $u$ be a positive solution of \eqref{1-1} .
Set $v=u^{-\frac{1}{\beta}}$, with  $\beta\neq0$  to be determined later and let $z=\left|\nabla v\right|^2$. Then
\begin{align}\label{5-1}
\Delta_mv=(\beta+1)(m-1)\frac{z^{\frac{m}{2}}}{v}+\frac{|\beta|^{2-m}}{\beta}v^\sigma+M|\beta|^{q-m}\beta v^sz^{\frac{q}{2}},
\end{align}
where
\begin{align}\label{sigma-s}
\begin{cases}
	\sigma=m-\beta(p-m+1)-1,\\
	s=(\beta+1)(m-q-1).
\end{cases}
\end{align}
By (\ref{5-1}), we obtain
\begin{align}\label{5-3}
z^{2-m}\left(\Delta_mv\right)^2&=(\beta+1)^2(m-1)^2\frac{z^2}{v^2}+\beta^{2(1-m)}v^{2\sigma}z^{2-m}
\nonumber\\
&\quad+M^2\beta^{2(q-m+1)}v^{2s}z^{q-m+2} +2M|\beta|^{q-2m+2}v^{\sigma+s}z^{\frac{q}{2}-m+2}
\nonumber\\&\quad
+2M|\beta|^{q-m}\beta (\beta+1)(m-1)v^{s-1}z^{\frac{q-m}{2}+2}
\nonumber\\
&\quad
+\frac{2|\beta|^{2-m}}{\beta}(\beta+1)(m-1)v^{\sigma-1}z^{2-\frac{m}{2}},
\end{align}
\begin{align}\label{5-4}
z^{1-\frac{m}{2}}\left<\nabla\Delta_mv,\nabla v\right>&=-(\beta+1)(m-1)\frac{z^2}{v^2}
+\frac{\sigma|\beta|^{2-m}}{\beta}v^{\sigma-1}z^{2-\frac{m}{2}}
\nonumber\\&\quad
+sM|\beta|^{q-m}\beta v^{s-1}z^{\frac{q-m}{2}+2}
\nonumber\\&\quad
+\frac{q}{2}M|\beta|^{q-m}\beta v^sz^{\frac{q-m}{2}}\left<\nabla z,\nabla v\right>
\nonumber\\&\quad
+\frac{m}{2}(\beta+1)(m-1)\frac{\left<\nabla z,\nabla v\right>}{v},
\end{align}
and
\begin{align}\label{5-5}
z^{-\frac{m}{2}}\Delta_mv\left<\nabla z,\nabla v\right>&=
(\beta+1)(m-1)\frac{\left<\nabla z,\nabla v\right>}{v}
\nonumber\\&\quad
+\frac{|\beta|^{2-m}}{\beta}v^\sigma z^{-\frac{m}{2}}\left<\nabla z,\nabla v\right>
\nonumber\\&\quad
+M|\beta|^{q-m}\beta v^sz^{\frac{q-m}{2}}\left<\nabla z,\nabla v\right>.
\end{align}
Substituting (\ref{5-3}), (\ref{5-4}) and (\ref{5-5}) into (\ref{L-1}), we derive
\begin{align}\label{5-6}
&\frac{1}{2}\mathscr A_v(z)
+\left(\frac{(\beta+1)(m-1)}{N}-1\right)(\beta+1)(m-1)\frac{z^2}{v^2}
\nonumber\\&\quad
+\left(\sigma+\frac{2(\beta+1)(m-1)}{N}\right)\frac{|\beta|^{2-m}}{\beta}v^{\sigma-1}z^{2-\frac{m}{2}}
\nonumber\\&\quad
+\left(s+\frac{2(\beta+1)(m-1)}{N}\right)M|\beta|^{q-m}\beta v^{s-1}z^{\frac{q-m}{2}+2}
\nonumber\\&\quad
+\frac{1}{N\beta^{2(m-1)}}v^{2\sigma}z^{2-m}
+\frac{M^2\beta^{2(q-m+1)}}{N}v^{2s}z^{q-m+2}
\nonumber\\&\quad
+\frac{2M|\beta|^{q-2m+2}}{N}v^{\sigma+s}z^{\frac{q}{2}-m+2}
\nonumber\\&\quad
-\frac{(N+2)(m-2)}{2N}\frac{|\beta|^{2-m}}{\beta}v^\sigma z^{-\frac{m}{2}}\left<\nabla z,\nabla v\right>
\nonumber\\&\quad
+\left(\frac{q}{2}-\frac{(N+2)(m-2)}{2N}\right)M|\beta|^{q-m}\beta v^sz^{\frac{q-m}{2}}\left<\nabla z,\nabla v\right>
\nonumber\\&\quad
+\left(\frac{m}{2}-\frac{(N+2)(m-2)}{2N}\right)(\beta+1)(m-1)\frac{\left<\nabla z,\nabla v\right>}{v}
\nonumber\\&\quad
+\frac{(2N+m-2)(m-2)}{4N}\frac{\left<\nabla z,\nabla v\right>^2}{z^2}
\nonumber\\&\quad
-\frac{m-2}{4}\frac{|\nabla z|^2}{z}
\leq0,\quad\mbox{on $\{z>0\}$}.
\end{align}
Afterwards, set $Y=v^\lambda z$ on $\{z>0\}$ for some parameter $\lambda$ to be determined later. In order to replace $\mathscr A_v(z)$ by $\mathscr A_v(Y)$,
we first calculate
\begin{align}\label{n-1}
-\Delta z&=\lambda v^{-\lambda-1}Y\Delta v-\lambda(\lambda+1)v^{-2\lambda-2}Y^2
\nonumber\\&\quad
+2\lambda v^{-\lambda-1}\left<\nabla v,\nabla Y\right>-v^{-\lambda}\Delta Y,
\end{align}
{where we have used that $v^{-2\lambda-2}Y|\nabla v|^2=v^{-\lambda-2}Y^2$.
Furthermore, reading the $m$-Laplacian as $\Delta_m v=\mbox{div}\bigl(z^{\frac{m}{2}-1}\nabla v\bigr)$},  we get
 \begin{align}
\Delta v
=z^{1-\frac{m}{2}}\Delta_mv-\frac{m-2}{2}\frac{\left<\nabla z,\nabla v\right>}{z}\nonumber.
\end{align}
Then using
\begin{equation}\label{scalar}
\left<\nabla z,\nabla v\right>=-\lambda v^{-2\lambda-1}Y^{2}+v^{-\lambda}\left<\nabla v,\nabla Y\right>,
\end{equation}
and \eqref{5-1}, we obtain
\begin{align}\label{5-7-jian}
\Delta v&=\left[\frac{\lambda(m-2)}{2}+(\beta+1)(m-1)\right]v^{-\lambda-1}Y
+\frac{|\beta|^{2-m}}{\beta}v^{\sigma-\lambda\left(1-\frac{m}{2}\right)}Y^{1-\frac{m}{2}}
\nonumber\\&\quad
+M|\beta|^{q-m}\beta v^{s-\lambda\left(\frac{q-m}{2}+1\right)}Y^{\frac{q-m}{2}+1}
-\frac{m-2}{2}\frac{\left<\nabla v,\nabla Y\right>}{Y}.
\end{align}
Replacing  \eqref{5-7-jian} into \eqref{n-1}, we obtain
\begin{align}\label{5-7}
-\Delta z&=\lambda\left[\lambda\left(\frac{m}{2}-2\right)+\beta(m-1)+m-2\right]v^{-2\lambda-2}Y^2
\nonumber\\&\quad
+\frac{\lambda|\beta|^{2-m}}{\beta}v^{\sigma-\lambda\left(2-\frac{m}{2}\right)-1}Y^{2-\frac{m}{2}}
\nonumber\\&\quad
+\lambda M|\beta|^{q-m}\beta v^{s-\lambda\left(\frac{q-m}{2}+2\right)-1}Y^{\frac{q-m}{2}+2}
\nonumber\\&\quad
+\lambda\left(3-\frac{m}{2}\right) v^{-\lambda-1}\left<\nabla v,\nabla Y\right>-v^{-\lambda}\Delta Y.
\end{align}
Next, we focus on $\frac{\left<D^2z\nabla v,\nabla v\right>}{z}$.
In view of (\ref{2-4-4}), we have
\begin{equation}\label{D^2z}
\left<D^2z\nabla v,\nabla v\right>=\left<\nabla\left<\nabla z,\nabla v\right>,\nabla v\right>-\frac{1}{2}|\nabla z|^2,
\end{equation}
and using \eqref{scalar} we get
\begin{align}
\left<\nabla\left<\nabla z,\nabla v\right>,\nabla v\right>&=
\lambda(2\lambda+1)v^{-3\lambda-2}Y^3
-3\lambda v^{-2\lambda-1}Y\left<\nabla v,\nabla Y\right>
\nonumber\\&\quad
+v^{-\lambda}\left<\nabla\left<\nabla v,\nabla Y\right>,\nabla v\right>,
\nonumber
\end{align}
 and, as in  \eqref{2-4} and using $\nabla z=2 D^2v \nabla v$,  we arrive to
$$
\begin{cases}
\left<\nabla\left<\nabla v,\nabla Y\right>,\nabla v\right>=
\left<D^2Y\nabla v,\nabla v\right>+\frac{1}{2}\left<\nabla z,\nabla Y\right>,\nonumber\\
\left<\nabla z,\nabla Y\right>=-\lambda v^{-\lambda-1}Y\left<\nabla v,\nabla Y\right>+v^{-\lambda}|\nabla Y|^2.
,\nonumber
\end{cases}
$$
Then,  by  \eqref{D^2z}, 
\begin{align}
\left<D^2z\nabla v,\nabla v\right>&=
\lambda(2\lambda+1)v^{-3\lambda-2}Y^3-\frac{7\lambda}{2}v^{-2\lambda-1}Y\left<\nabla v,\nabla Y\right>
\nonumber\\&\quad
+v^{-\lambda}\left<D^2Y\nabla v,\nabla v\right>
+\frac{1}{2}v^{-2\lambda}|\nabla Y|^2
-\frac{1}{2}|\nabla z|^2.\nonumber
\end{align}
Thus
\begin{align}\label{5-8}
\frac{\left<D^2z\nabla v,\nabla v\right>}{z}&=
\lambda(2\lambda+1)v^{-2\lambda-2}Y^2-\frac{7\lambda}{2} v^{-\lambda-1}\left<\nabla v,\nabla Y\right>
\nonumber\\&\quad
+\frac{\left<D^2Y\nabla v,\nabla v\right>}{Y}
+\frac{1}{2}v^{-\lambda}\frac{|\nabla Y|^2}{Y}
-\frac{1}{2}\frac{|\nabla z|^2}{z}.
\end{align}
Combining (\ref{5-7}) and (\ref{5-8}), we derive
\begin{align}\label{5-9}
\mathscr A_v(z)&=-\Delta z-(m-2)\frac{\left<D^2z\nabla v,\nabla v\right>}{z}
\nonumber\\
&=v^{-\lambda}\mathscr A_v(Y)+\lambda\left[\lambda\left(2-\frac{3m}{2}\right)+\beta(m-1)\right]v^{-2\lambda-2}Y^2
\nonumber\\&\quad
+\frac{\lambda|\beta|^{2-m}}{\beta}v^{\sigma-\lambda\left(2-\frac{m}{2}\right)-1}Y^{2-\frac{m}{2}}
\nonumber\\&\quad
+\lambda M|\beta|^{q-m}\beta v^{s-\lambda\left(\frac{q-m}{2}+2\right)-1}Y^{\frac{q-m}{2}+2}
\nonumber\\&\quad
+\lambda(3m-4) v^{-\lambda-1}\left<\nabla v,\nabla Y\right>
\nonumber\\&\quad
-\frac{m-2}{2}v^{-\lambda}\frac{|\nabla Y|^2}{Y}
+\frac{m-2}{2}\frac{|\nabla z|^2}{z}.
\end{align}
Replacing  into \eqref{5-6}, the following expressions
$$\begin{aligned}
\frac{z^2}{v^2}&=v^{-2\lambda-2}Y^2,\\
v^{\sigma-1}z^{2-\frac{m}{2}}&=v^{\sigma-\lambda\left(2-\frac{m}{2}\right)-1}Y^{2-\frac{m}{2}},\\
v^{s-1}z^{\frac{q-m}{2}+2}&=v^{s-\lambda\left(\frac{q-m}{2}+2\right)-1}Y^{\frac{q-m}{2}+2},\\
v^{2\sigma}z^{2-m}&=v^{2\sigma-\lambda(2-m)}Y^{2-m},\\
v^{2s}z^{q-m+2}&=v^{2s-\lambda(q-m+2)}Y^{q-m+2},\\
v^{\sigma+s}z^{\frac{q}{2}-m+2}&=v^{\sigma+s-\lambda\left(\frac{q}{2}-m+2\right)}Y^{\frac{q}{2}-m+2},\\
v^\sigma z^{-\frac{m}{2}}\left<\nabla z,\nabla v\right>&=-\lambda v^{\sigma-\lambda\left(2-\frac{m}{2}\right)-1}Y^{2-\frac{m}{2}}
+v^{\sigma-\lambda\left(1-\frac{m}{2}\right)}Y^{-\frac{m}{2}}\left<\nabla v,\nabla Y\right>,
\\
v^sz^{\frac{q-m}{2}}\left<\nabla z,\nabla v\right>&=-\lambda v^{s-\lambda\left(\frac{q-m}{2}+2\right)-1}Y^{\frac{q-m}{2}+2}
+v^{s-\lambda\left(\frac{q-m}{2}+1\right)}Y^{\frac{q-m}{2}}\left<\nabla v,\nabla Y\right>,
\\
\frac{\left<\nabla z,\nabla v\right>}{v}&=-\lambda v^{-2\lambda-2}Y^{2}+v^{-\lambda-1}\left<\nabla v,\nabla Y\right>,\\
\frac{\left<\nabla z,\nabla v\right>^2}{z^2}&=\lambda^2v^{-2\lambda-2}Y^2-2\lambda v^{-\lambda-1}\left<\nabla v,\nabla Y\right>+\frac{\left<\nabla v,\nabla Y\right>^2}{Y^2},
\end{aligned} $$
 we get an estimate from above  for $\mathscr A_v(z)$, precisely
 \begin{align}\label{A_estim}
 	\mathscr A_v(z)&\le\,\biggr\{2(\beta+1)(m-1)\biggl[\lambda\left(\frac{m}{2}-\frac{(N+2)(m-2)}{2N}\right)-\left(\frac{(\beta+1)(m-1)}{N}-1\right)\biggr]\nonumber\\
 	&\qquad -\lambda^2\biggl(m-2+\frac{(m-2)^2}{2N}\biggr)\biggr\}v^{-2\lambda-2}Y^2
 	\nonumber\\
 	&\quad-\biggl[2\sigma +\frac{4(\beta+1)(m-1)}N+\lambda \frac{(N+2)(m-2)}{N}\biggr] \frac{|\beta|^{2-m}}{\beta} v^{\sigma-\lambda\left(2-\frac{m}{2}\right)-1}Y^{2-\frac{m}{2}}
 	\nonumber\\&\quad
 	-\biggl[2s+\frac{4(\beta+1)(m-1)}N-\lambda \biggl(q-\frac{(N+2)(m-2)}{N}\biggr)\biggr]\cdot\nonumber\\
	&\qquad\qquad\qquad\qquad\qquad\qquad\qquad\qquad\cdot M|\beta|^{q-m}\beta v^{s-\lambda\left(\frac{q-m}{2}+2\right)-1}Y^{\frac{q-m}{2}+2}
 	\nonumber\\&\quad
 	-\frac{2}{N\beta^{2(m-1)}}v^{2\sigma-\lambda(2-m)}Y^{2-m}
 	-2\frac{M^2\beta^{2(q-m+1)}}{N}v^{2s-\lambda(q-m+2)}Y^{q-m+2}
 	\nonumber\\&\quad
 	-\frac{4M}{N}|\beta|^{q-2m+2}v^{\sigma+s-\lambda\left(\frac{q}{2}-m+2\right)}Y^{\frac{q}{2}-m+2}
 	\nonumber\\
 	&\quad
 	-\biggl[ 2(\beta+1)(m-1)\biggl(1-\frac{m-2}N\biggr)-\!\lambda(m-2)\biggl(2+\frac{m-2}N\biggr)\biggr]v^{-\lambda-1}\left<\nabla v,\nabla Y\right>
 	\nonumber\\&\quad
 	+\frac{(N+2)(m-2)}{N}\frac{|\beta|^{2-m}}{\beta}v^{\sigma-\lambda\left(1-\frac{m}{2}\right)}Y^{-\frac{m}{2}}\left<\nabla v,\nabla Y\right>
 	\nonumber\\&\quad
 	-\left(q-\frac{(N+2)(m-2)}{N}\right)M|\beta|^{q-m}\beta v^{s-\lambda\left(\frac{q-m}{2}+1\right)}Y^{\frac{q-m}{2}}\left<\nabla v,\nabla Y\right>
 	\nonumber\\&\quad
 	-\frac{(2N+m-2)(m-2)}{2N}\frac{\left<\nabla v,\nabla Y\right>^2}{Y^2}
 	\nonumber\\&\quad
 	+\frac{m-2}{2}\frac{|\nabla z|^2}{z},\qquad\mbox{for $z>0$}.
 \end{align}
Replacing
 \eqref{5-9} in \eqref{A_estim} we deduce that for some positive constant $c_6=c_6(N,m,q,\beta,\lambda)$, the following holds
\begin{align}\label{5-10}
&\quad v^{-\lambda}\mathscr A_v(Y)+
L_1v^{-2\lambda-2}Y^2+L_2v^{\sigma-\lambda\left(2-\frac{m}{2}\right)-1}Y^{2-\frac{m}{2}}
\nonumber\\&\quad
+L_3v^{s-\lambda\left(\frac{q-m}{2}+2\right)-1}Y^{\frac{q-m}{2}+2}
+L_4v^{2\sigma-\lambda(2-m)}Y^{2-m}
\nonumber\\&\quad
+L_5v^{2s-\lambda(q-m+2)}Y^{q-m+2}
+L_6v^{\sigma+s-\lambda\left(\frac{q}{2}-m+2\right)}Y^{\frac{q}{2}-m+2}
\nonumber\\
&\leq
c_6\left\{\left(v^{-\lambda-1}
+v^{\sigma-\lambda\left(1-\frac{m}{2}\right)}Y^{-\frac{m}{2}}\right)\left|\left<\nabla v,\nabla Y\right>\right|+\frac{\left<\nabla v,\nabla Y\right>^2}{Y^2}\right.
\nonumber\\&\left.\quad
+v^{-\lambda}\frac{|\nabla Y|^2}{Y}\right\}+L_7v^{s-\lambda\left(\frac{q-m}{2}+1\right)}Y^{\frac{q-m}{2}}\left|\left<\nabla v,\nabla Y\right>\right|,
\end{align}
where
\begin{align}
L_1&=\lambda^2\left(\frac{(m-2)^2}{2N}-\frac{m}{2}\right)-\lambda(m-1)\left(\beta+2-\frac{2(\beta+1)(m-2)}{N}\right)
\nonumber\\&\quad
+2(\beta+1)(m-1)\left(\frac{(\beta+1)(m-1)}{N}-1\right),
\nonumber\\
L_2&=\frac{|\beta|^{2-m}}{\beta}\left\{\lambda\left(m-1+\frac{2(m-2)}{N}\right)+\frac{4(\beta+1)(m-1)}{N}+2\sigma\right\},
\nonumber\\
L_3&=M|\beta|^{q-m}\beta\left\{\lambda\left(m-q-1+\frac{2(m-2)}{N}\right)+\frac{4(\beta+1)(m-1)}{N}+2s\right\},
\nonumber\\
L_4&=\frac{2}{N\beta^{2(m-1)}},\quad L_5=\frac{2M^2\beta^{2(q-m+1)}}{N},\quad L_6=\frac{4M|\beta|^{q-2m+2}}{N},
\nonumber\\
L_7&=\left(q+\frac{(N+2)|m-2|}{N}\right)M|\beta|^{q-m+1}.
\nonumber
\end{align}
In particular, it results that $L_4,\,L_5,\,L_6,\,L_7>0$.

Multiplying (\ref{5-10}) by $v^\lambda$ yields
\begin{align}\label{5-11}
&\quad \mathscr A_v(Y)+
L_1v^{-\lambda-2}Y^2+L_2v^{\sigma-\lambda\left(1-\frac{m}{2}\right)-1}Y^{2-\frac{m}{2}}
\nonumber\\&\quad
+L_3v^{s-\lambda\left(\frac{q-m}{2}+1\right)-1}Y^{\frac{q-m}{2}+2}
+L_4v^{2\sigma+\lambda(m-1)}Y^{2-m}
\nonumber\\&\quad
+L_5v^{2s-\lambda(q-m+1)}Y^{q-m+2}
+L_6v^{\sigma+s-\lambda\left(\frac{q}{2}-m+1\right)}Y^{\frac{q}{2}-m+2}
\nonumber\\&
\leq c_6\left\{\left(v^{-1}
+v^{\sigma+\frac{m\lambda}{2}}Y^{-\frac{m}{2}}\right)\left|\left<\nabla v,\nabla Y\right>\right|+v^\lambda \frac{\left<\nabla v,\nabla Y\right>^2}{Y^2}\right.
\nonumber\\&\left.\quad
+\frac{|\nabla Y|^2}{Y}\right\}+L_7v^{s-\lambda\frac{q-m}{2}}Y^{\frac{q-m}{2}}\left|\left<\nabla v,\nabla Y\right>\right|.
\end{align}
Now we estimate each term in the right-hand side of (\ref{5-11}).
For any $\varepsilon>0$, using that $|\nabla v|^2=v^{-\lambda}Y$, we have
$$c_6\frac{\left|\left<\nabla v,\nabla Y\right>\right|}{v} \le v^{-\frac\lambda 2-1}\sqrt{Y}|\nabla Y|\le\varepsilon v^{-\lambda-2}Y^2+\frac{c_6^2}{ 4\varepsilon}\frac{|\nabla Y|^2}{Y}.$$
and
$$\begin{aligned}
	c_6v^{\sigma+\frac{m\lambda}{2}}Y^{-\frac{m}{2}}\left|\left<\nabla v,\nabla Y\right>\right|&\le c_6 v^{\sigma+\frac \lambda 2 (m-1)}Y^{\frac{2-m}2-\frac12 }|\nabla Y|
\\&\leq  \varepsilon v^{2\sigma+\lambda(m-1)}Y^{2-m}+\frac{c_6^2}{4\varepsilon}\frac{|\nabla Y|^2}{Y}.\end{aligned}$$
Similarly, being $L_5$ positive, we get 
$$L_7v^{s-\lambda\frac{q-m}{2}}Y^{\frac{q-m}{2}}\left|\left<\nabla v,\nabla Y\right>\right|\leq \frac{L_5}{2} v^{2s-\lambda(q-m+1)}Y^{q-m+2}+\frac{L_7^2}{2L_5}\frac{|\nabla Y|^2}{Y}.$$
Noting also that
$$v^\lambda \frac{\left<\nabla v,\nabla Y\right>^2}{Y^2}\leq \frac{|\nabla Y|^2}{Y}.$$
Hence, we obtain from (\ref{5-11}) that
\begin{align}\label{5-12}
	\mathscr A_v(Y)+H_1+H_2\leq c_7\frac{|\nabla Y|^2}{Y},
\end{align}
where $c_7=c_7(N,m,q,\beta,\lambda)>0$, and
\begin{align}\label{5-13}
H_1&:=(L_1-\varepsilon)v^{-\lambda-2}Y^2+L_2v^{\sigma-\lambda\left(1-\frac{m}{2}\right)-1}Y^{2-\frac{m}{2}}
\nonumber\\&\quad +  (L_4-\varepsilon)  v^{2\sigma+\lambda(m-1)}Y^{2-m}\nonumber\\
& =v^{-\lambda-2}Y^2\biggl[ L_1-\varepsilon+L_2v^{\sigma+\lambda\frac{m}{2}+1}Y^{-\frac{m}{2}}+ (L_4-\varepsilon)  v^{2\sigma+\lambda m+2}Y^{-m}\biggr] ,
\end{align}
and
\begin{align}\label{5-14}
H_2&:=L_3v^{s-\lambda\left(\frac{q-m}{2}+1\right)-1}Y^{\frac{q-m}{2}+2}
+\frac{L_5}{2}v^{2s-\lambda(q-m+1)}Y^{q-m+2}
\nonumber\\&\quad
+L_6v^{\sigma+s-\lambda\left(\frac{q}{2}-m+1\right)}Y^{\frac{q}{2}-m+2}.
\end{align}
Now, fix
$$\lambda<-2,\qquad \beta>0, \qquad 2(\beta+1)+\lambda>0.$$
By this choice,  we immediately see that the positivity of $H_2$ is ensured,  indeed the second and the
third terms of $H_2$ are positive, being $L_5,\,L_6>0$,  it remains to prove that $L_3>0$. This latter follows by the positivity of
$$L_3'=\lambda\left(m-q-1+\frac{2(m-2)}{N}\right)+\frac{4(\beta+1)(m-1)}{N}+2s.$$ 
 Since, $s=(\beta+1)(m-q-1)$,  by \eqref{sigma-s},  and $m-1<q<\frac{(N+2)(m-1)}{N}$,  by assumption,  then we obtain
$$L_3'=\frac{2\lambda(m-2)}{N}+\frac{4(\beta+1)(m-1)}{N}-(q-m+1)\left[2(\beta+1)+\lambda\right]>-\frac{2\lambda}{N}>0.$$
To estimate the term $H_1$, we consider the following  trinomial
\begin{align}
T_\varepsilon(t)=  (L_4-\varepsilon) t^2+L_2t+L_1-\varepsilon.\nonumber
\end{align}
If its discriminant is strictly negative, then it is possible to find $\gamma $  small enough so that the discriminant of  $(L_4-\varepsilon-\gamma) t^2+L_2t+L_1-\varepsilon-\gamma$ still remains strictly negative, in turn we can conclude that
there exists $\gamma=\gamma(N,m,p,q,\beta,\lambda,\varepsilon)>0$ such that $T_{\varepsilon}(t)\geq\gamma\left(t^2+1\right)$,
 and hence
\begin{align}
H_1&=v^{-\lambda-2}Y^2T_\varepsilon\left(v^{\sigma+\frac{m\lambda}{2}+1}Y^{-\frac{m}{2}}\right)
\nonumber\\&
\geq\gamma\left(v^{-\lambda-2}Y^2+v^{2\sigma+\lambda(m-1)}Y^{2-m}\right)\nonumber
\end{align}
Since $\lambda<-2$, we can define
\begin{align}
S=\frac{2\sigma+\lambda(m-1)}{\lambda+2}=m-1+\frac{p-m+1}{d}.\nonumber
\end{align}
where  $d:=-\frac{\lambda+2}{2\beta}>0$,  by the choice of $\lambda$ and $\beta$, so that  $S>m-1$.

Since  $\frac{2S-m+2}{S+1}>1$, we have
\begin{align}
Y^{\frac{2S-m+2}{S+1}}&=\left(v^{-\lambda-2}Y^2\right)^{\frac{S}{S+1}}\left(v^{(\lambda+2)S}
Y^{2-m}\right)^{\frac{1}{S+1}}
\nonumber\\&
\leq v^{-\lambda-2}Y^2+v^{(\lambda+2)S}Y^{2-m}
\nonumber\\&
=v^{-\lambda-2}Y^2+v^{2\sigma+\lambda(m-1)}Y^{2-m}.\nonumber
\end{align}
Therefore,
\begin{align}\label{5-15}
H_1\geq\gamma Y^{\frac{2S-m+2}{S+1}}.
\end{align}

Combining with \eqref{5-12}, \eqref{5-13}, \eqref{5-14} and \eqref{5-15}, we arrive
\begin{align}\label{5-16}
\mathscr A_v(Y)+\gamma Y^{\frac{2S-m+2}{S+1}}\leq c_7\frac{\left|\nabla Y\right|^2}{Y}.
\end{align}
By Lemma \ref{lem1}, we obtain
$$Y(x)\leq c_8\left({\rm dist}\left(x,\partial\Omega\right)\right)^{-\frac{2(S+1)}{S-m+1}}=c_8\left({\rm dist}\left(x,\partial\Omega\right)\right)^{-\frac{2\sigma+m\lambda+2}{\sigma-m+1}},$$
where $c_8=c_8(S,m,\gamma,c_7)>0$. It follows that
\begin{align}\label{5-17}
\left|\nabla u^{d}(x)\right|\leq c_8'\left({\rm dist}\left(x,\partial\Omega\right)\right)^{-\frac{2\sigma+m\lambda+2}{2(\sigma-m+1)}}
= c_8'\left({\rm dist}\left(x,\partial\Omega\right)\right)^{-1-\frac{md}{p-m+1}},
\end{align}
where $c_8'=c_8'(m,\lambda,\beta)>0$,
which is exactly (\ref{1-8}).
The nonexistence of any positive solution of (\ref{1-1}) in $\mathbb R^N$ follows consequently.

It remains to prove that the discriminant of the trinomial $T_\varepsilon(t)$ is negative. Since the discriminant is a polynomial of its coefficients.
Hence it suffices to prove that the discriminant of $T_0(t)$ is  strictly negative to deduce the same property holds for $T_\varepsilon(t)$ for small enough $\varepsilon$. Noting that
\begin{align}
T_0(t)=L_4t^2+L_2t+L_1,\nonumber
\end{align}
and its discriminant $D=L_2^2-4L_1L_4$ satisfies
\begin{align}
D&=|\beta|^{2(1-m)}\left\{\left[\lambda\left(m-1+\frac{2(m-2)}{N}\right)+\frac{4(\beta+1)(m-1)}{N}+2\sigma\right]^2\right.
\nonumber\\&\quad
-\frac{8\lambda^2}{N}\left(\frac{(m-2)^2}{2N}-\frac{m}{2}\right)
+\frac{8\lambda(m-1)}{N}\left(\beta+2-\frac{2(\beta+1)(m-2)}{N}\right)
\nonumber\\&\quad\left.
-\frac{16}{N}(\beta+1)(m-1)\left(\frac{(\beta+1)(m-1)}{N}-1\right)\right\}.
\nonumber
\end{align}
Using $\beta+1=\frac{2p+\lambda(m-1)-(\lambda+2)S}{2(p-m+1)}$ and $\sigma=\frac{(\lambda+2)S-\lambda(m-1)}{2}$, we further compute
$$ \begin{aligned}\lambda&\left(m-1+\frac{2(m-2)}{N}\right)+\frac{4(\beta+1)(m-1)}{N}+2\sigma=\frac1{N(p-m+1)}
\\&\quad\times\biggl\{4p(m-1)+2\lambda\bigl[m-1+p(m-2)\bigr] +(\lambda+2)S\bigl[N(p-m+1)-2(m-1)\bigr]\biggr\},
\end{aligned}$$
$$\beta+2-\frac{2(\beta+1)(m-2)}{N}=\frac{2p+\lambda(m-1)-(\lambda+2)S}{2N(p-m+1)}[N-2(m-2)]+1,$$
and
$$\begin{aligned}(\beta+1)(m-1)&\left(\frac{(\beta+1)(m-1)}{N}-1\right)=\frac{(m-1)[2p+\lambda(m-1)-(\lambda+2)S]}{4N(p-m+1)^2}\\&\times \bigl\{(m-1)[2p+\lambda(m-1)-(\lambda+2)S]-2N(p-m+1)\bigr\}.
\end{aligned}$$
Thus
\begin{align}
D&=\frac{\beta^{2(1-m)}}{N(p-m+1)}\left\{(\lambda+2)^2\left[N(p-m+1)-4(m-1)\right]S^2\right.
\nonumber\\&\quad
+4(\lambda+2)\left[\lambda p(m-2)+2(m-1)(p-1)\right]S
\nonumber\\&\quad\left.
+4\lambda^2(p-m+1)+4(\lambda+2)^2p(m-1)\right\}.
\nonumber
\end{align}
 Since $\lambda+2\neq0$,  we set $\ell=\frac{\lambda}{\lambda+2}$. By the choice of $\lambda$ it follows $\ell>1$. In turn, using also that $1/(\lambda+2)=(1-\ell)/2<0$, we arrive to
\begin{align}
D&=\frac{(\lambda+2)^2\beta^{2(1-m)}}{N(p-m+1)}\left\{\left[N(p-m+1)-4(m-1)\right]S^2-4(p-m+1)\ell S\right.
\nonumber\\&\quad\left.
+4(m-1)(p-1)S+4(p-m+1)\ell^2+4p(m-1)\right\},
\nonumber
\end{align}
which is equivalent to
\begin{align}
D&=\frac{(\lambda+2)^2\beta^{2(1-m)}}{N(p-m+1)}\left\{4(p-m+1)\left(\ell-\frac{S}{2}\right)^2+D_1(S)\right\},
\nonumber
\end{align}
where $$D_1(S):=\left[(N-1)(p-m+1)-4(m-1)\right]S^2+4(m-1)(p-1)S+4p(m-1).$$
Fix $\ell=\frac{S}{2}$, hence $\beta=\frac{\lambda(m-3)+2(m-1)}{2(p-m+1)}$.
As the coefficient of $S^2$ in $D_1(S)$ is negative if $p<\frac{(N+3)(m-1)}{N-1}$, we can choose $S$ large enough, namely  $\lambda<-2$ such that  $|\lambda+2|$ is small enough,
to reach $D_1(S)<0$.  In particular,  condition $2(\beta+1)+\lambda>0$ holds true for $\lambda\to-2^-$ being equivalent to  $(\lambda+2)p-2\lambda>0$.
Consequently  $D<0$, concluding the proof of  the positivity of $T_{\varepsilon}$. 
\end{proof}

\end{document}